\documentclass[11pt,reqno]{amsart}

\usepackage{amsmath,amssymb,amsthm,mathtools}
\usepackage{enumitem}
\usepackage{lmodern}
\usepackage{microtype}
\usepackage{xcolor}
\usepackage[hidelinks]{hyperref}
\hypersetup{
 pdftitle={Rigidity and Classification of Manifolds with Maximal Boundary Static Kernels},
 pdfauthor={Hongyi Sheng and Caihong Yi}
}

\numberwithin{equation}{section}

\newtheorem{theorem}{Theorem}[section]
\newtheorem{proposition}[theorem]{Proposition}
\newtheorem{lemma}[theorem]{Lemma}
\newtheorem{corollary}[theorem]{Corollary}
\theoremstyle{definition}
\newtheorem{definition}[theorem]{Definition}
\newtheorem{example}[theorem]{Example}
\theoremstyle{remark}
\newtheorem{remark}[theorem]{Remark}

\newcommand{\R}{\mathbb R}
\newcommand{\Sph}{\mathbb S}
\newcommand{\Hyp}{\mathbb H}
\newcommand{\Mspace}{\mathbb M}
\newcommand{\Ric}{\operatorname{Ric}}
\newcommand{\Hess}{\operatorname{Hess}}
\newcommand{\tr}{\operatorname{tr}}
\newcommand{\Span}{\operatorname{span}}

\newcommand{\Lie}{\mathcal L}

\newcommand{\inner}[2]{\left\langle #1,#2\right\rangle}

\title[Manifolds with Maximal Boundary Static Kernels]
{Rigidity and Classification of Manifolds with Maximal Boundary Static Kernels}

\author[H. Sheng]{Hongyi Sheng}
\address{Institute for Theoretical Sciences, Westlake Institute for Advanced Study,
Westlake University, Hangzhou, Zhejiang Province 310030, China}
\email{shenghongyi@westlake.edu.cn}

\author[C. Yi]{Caihong Yi}
\address{School of Mathematics, Hangzhou Normal University, Hangzhou 311121, China}
\email{caihongyi@hznu.edu.cn}



\begin{document}
\raggedbottom

\begin{abstract}
Boundary static potentials arise from the adjoint equations associated with
simultaneous variations of scalar curvature and boundary mean
curvature. On a connected manifold of dimension
$n\ge2$ with nonempty smooth boundary, the space of boundary static
potentials has dimension at most $n$. We classify all complete
manifolds attaining this bound. For $n\ge3$, every such manifold is
isometric to a domain in a simply connected space form and has at
most two boundary components. In dimension two, the classification
also includes coverings of annuli and exteriors, as well as
spherical M\"obius bands. Finally, for $n\ge3$, the dimension hypothesis is sharp for
constant-curvature rigidity, and for $n\ge2$, it is sharp for
the bound on the number of boundary components.
\end{abstract}

\maketitle

\section{Introduction}
The deformation theory of scalar curvature plays a central role in
geometric analysis and general relativity, particularly through the
interplay between the scalar curvature map and static metrics.
In the closed and localized settings, the absence of static
potentials yields local surjectivity of the scalar curvature map
at $g$ \cite{FischerMarsden1975,Corvino2000}.
On the other hand, the existence of static potentials is closely
related to rigidity. In particular, Qing and Yuan proved a local
scalar curvature rigidity theorem for sufficiently small geodesic
balls centered at points where a static potential is nonzero
\cite{QingYuan2016}.
More broadly, the classification of static manifolds is a
longstanding problem, with classical results due to Kobayashi
and Lafontaine \cite{Kobayashi1982,Lafontaine1983}.

On manifolds with boundary, simultaneous variations of scalar
curvature and boundary mean curvature lead to an additional
geometric boundary condition
\cite{CruzVitorio2019,ShengLocalized2024}.
The resulting adjoint system consists of the interior static
equation together with a boundary condition involving the
second fundamental form. The existence of a nontrivial solution
therefore imposes restrictions on both the interior geometry
and the geometry of the boundary.
Understanding these restrictions leads naturally to the rigidity
and classification of manifolds admitting such potentials.
To make this discussion precise, we adopt the terminology
of \cite{ShengZhao2026}, allowing dimension two and imposing
completeness only when explicitly stated.

\begin{definition}
Let $(M^n,g)$, $n\ge2$, be a connected smooth Riemannian manifold
with nonempty smooth boundary.
A nontrivial function $u\in C^\infty(M)$ is called a
\emph{static potential} if it satisfies
\begin{equation}\label{static int}
 -(\Delta u)g+\operatorname{Hess}u-u\operatorname{Ric}=0
 \quad\text{in }M.
\end{equation}
A metric $g$ admitting such a potential is said to be
\emph{static}.
A static potential $u$ is called a \emph{boundary static potential}
(or \emph{b-static potential}) if
\begin{equation}\label{static bdry}
 u_\nu\hat g=uh
 \quad\text{on }\partial M,
\end{equation}
where $u_\nu$ denotes the derivative in the direction of the
outward unit normal $\nu$, $\hat g$ is the induced metric on
$\partial M$, and $h$ is the second fundamental form.
A manifold admitting such a potential is called \emph{b-static}.
\end{definition}

To relate the preceding equations to curvature deformation,
define the operator $L^*$ by
\[
 L^*u:=-(\Delta u)g+\operatorname{Hess}u-u\operatorname{Ric}.
\]
This is the operator in \eqref{static int} and the formal
adjoint of the linearization of scalar curvature.
Incorporating the boundary condition, we set
\[
 \Phi^*u:=(L^*u,u_\nu\hat g-uh).
\]
With this notation, the static potentials and b-static
potentials are precisely the nonzero elements of $\ker L^*$
and $\ker\Phi^*$, respectively.
Both kernels are taken in $C^\infty(M)$, and
$\ker\Phi^*\subseteq\ker L^*$.
We suppress subscripts for the fixed ambient metric $g$,
retaining them only when needed to make the metric dependence
explicit.
Recall that, in the compact and localized settings,
$\ker\Phi_g^*=\{0\}$
implies local surjectivity of $g\mapsto(R_g,H_g)$
\cite[Theorem~1.1]{ShengLocalized2024}.

When $\ker\Phi^*\ne\{0\}$, $M$ has constant scalar curvature $R$
and totally umbilic boundary, with constant mean curvature $H$
on each connected boundary component.
Moreover, $\operatorname{Ric}(X,\nu)=0$ for every
$X\in T(\partial M)$; see
\cite{ShengLocalized2024,Sheng2025} and the related rigidity
results in \cite{CruzNunes2023}.
The basic structural properties above follow from the existence
of a single b-static potential.
It is therefore natural to ask what
further geometric restrictions arise when several linearly
independent b-static potentials exist.

To address this question, we first recall the dimension bounds
for these solution spaces.
A static potential is determined by its value and differential
at one point, giving
\[
 \dim\ker L^*\le n+1;
\]
see \cite[Corollary~2.4]{Corvino2000}.
For a b-static potential, the boundary condition determines
its normal derivative from its boundary value. So
its value and tangential differential at a boundary point
determine the entire solution, and
\begin{equation*}
 \dim\ker\Phi^*\le n;
\end{equation*}
see \cite[Proposition~3.4]{ShengLocalized2024}.
Neither bound requires completeness.

The dimension of these spaces is also related to the symmetries
of the metric. Indeed, given $u,v\in\ker L^*$, the vector field
\begin{equation}\label{eq:intro-Xuv}
 X_{u,v}=uD v-vD u
\end{equation}
is Killing; see \cite{Lafontaine2009, CerneaGuan2015}.
If $u,v\in\ker\Phi^*$, the boundary condition \eqref{static bdry} also makes
$X_{u,v}$ tangent to $\partial M$.
Thus pairs of b-static potentials generate infinitesimal
isometries that preserve the boundary, suggesting stronger
rigidity when the kernel has maximal dimension.

Motivated by this relationship between kernel dimension and
symmetry, we study the maximal case
\[
 \dim\ker\Phi^*=n.
\]
Under this hypothesis, we classify all complete connected
manifolds with nonempty smooth boundary in every dimension
$n\ge2$.

For $n\ge3$, we first prove that maximality forces constant
sectional curvature.
This \emph{local conclusion}, however, does not determine the
topology of a manifold with boundary, the number of its
boundary components, or their arrangement.
Completeness allows us to resolve these \emph{global questions}
without assuming an ambient embedding, simple connectedness,
boundary topology, or an asymptotic model.
For surfaces, a separate analysis is needed to determine the possible global identifications.

To place the local curvature statement in context, recall that
Guan and Cern\'ea proved that a closed connected orientable
manifold of dimension $n\ge3$ admitting $n$ linearly independent
static potentials must be round \cite{GuanCernea2016}.
A related dimension gap follows from the warped product structure
theorem of He, Petersen and Wylie
\cite{HePetersenWylie2015}.
Our interior rigidity theorem gives a direct pointwise proof that
$n$ independent static potentials force constant sectional curvature
when $n\ge3$, without completeness or simple connectedness
(Theorem~\ref{thm:interior-maximality}).
Since $\ker\Phi^*\subseteq\ker L^*$, it applies under our
maximality assumption. Combining this result with the boundary
equation and the boundary one-jet isomorphism gives the following
description of the boundary geometry.
Here, $\Mspace_\kappa^n$ denotes the simply connected complete
space form of dimension $n$ and sectional curvature $\kappa$.

\begin{theorem}[Curvature rigidity and boundary geometry]\label{thm:local-main}
Let $(M^n,g)$, $n\ge3$, be a connected Riemannian
manifold with nonempty smooth boundary. If $\dim\ker\Phi^*=n$, then
$g$ has constant sectional curvature $\kappa=\frac{R}{n(n-1)}$. For every
boundary component $\Sigma$, the following hold:
\begin{enumerate}[label=\textup{(\roman*)}]
\item $\Sigma$ is totally umbilic,
$\lambda=\frac{H}{n-1}$ is constant, and the induced metric has
constant sectional curvature $k_\Sigma=\kappa+\lambda^2$.
\item Restriction to $\Sigma$ identifies $\ker\Phi^*$
with the full solution space of
\[
 \Hess_{\hat g}\phi=-k_\Sigma\phi\hat g
 \quad\text{on }\Sigma.
\]
\item If $\Sigma$ is complete, then
$\Sigma\cong\Mspace_{k_\Sigma}^{n-1}$. In particular, if
$\Sigma$ is compact, it is a round sphere.
\end{enumerate}
\end{theorem}

We now turn to the global classification, where completeness
becomes essential. Indeed,
Example~\ref{ex:incomplete-many-components} gives incomplete
manifolds with $\dim\ker\Phi^*=n$ and arbitrarily many boundary
components. We begin with the compact case, in which
completeness is automatic.

\begin{theorem}[Compact classification]\label{thm:main}
Let $(M^n,g)$, $n\ge3$, be a compact connected
Riemannian manifold with nonempty smooth boundary. If
$\dim\ker\Phi^*=n$, then $g$ has constant sectional curvature
$\kappa=\frac{R}{n(n-1)}$, and exactly one of the following holds:
\begin{enumerate}[label=\textup{(\roman*)}]
\item $\partial M$ is connected, and $M$
is isometric to a closed geodesic ball in $\Mspace_\kappa^n$;
\item $\partial M$ has two components, and
$M$ is isometric to the closed region between two concentric
geodesic spheres in $\Mspace_\kappa^n$.
\end{enumerate}
\end{theorem}

We next consider complete noncompact manifolds with at least
one compact boundary component.

\begin{theorem}[Noncompact case with a compact boundary component]\label{thm:noncompact-main}
Let $(M^n,g)$, $n\ge3$, be a complete connected
noncompact Riemannian manifold with smooth boundary. Suppose that
$\partial M$ has a compact component and $\dim\ker\Phi^*=n$.
Then $M$ is isometric to the closed exterior
\[
 \Mspace_\kappa^n\setminus B_{r_0}(o)
 =\{x\in\Mspace_\kappa^n:d(o,x)\ge r_0\},\qquad
 \kappa\le0,\quad r_0>0,
\]
of a geodesic ball in Euclidean or hyperbolic space. In particular,
$\partial M$ is connected and compact, and $M$ has exactly
one end.
\end{theorem}

It remains to consider manifolds whose boundary components are
all noncompact. By Theorem~\ref{thm:noncompact-main}, compact
and noncompact boundary components cannot coexist under
completeness and maximality. It therefore suffices to assume
that one component is noncompact.

\begin{theorem}[Noncompact boundary classification]\label{thm:noncompact-boundary-main}
Let $(M^n,g)$, $n\ge3$, be a complete connected
Riemannian manifold with smooth boundary. Suppose that
$\partial M$ has a noncompact component and $\dim\ker\Phi^*=n$.
Then $M$ is isometric to one of the following warped products:
\begin{enumerate}[label=\textup{(\roman*)}]
\item $I\times\R^{n-1}$ with $g=dt^2+g_{\R^{n-1}}$;
\item $I\times\R^{n-1}$ with $g=dt^2+e^{2at}g_{\R^{n-1}}$,
where $a>0$ is a constant;
\item $I\times\Hyp^{n-1}_{-a^2}$ with
$g=dt^2+\cosh^2(at)g_{\Hyp^{n-1}_{-a^2}}$, where $a>0$ is a constant.
\end{enumerate}
Here $\Hyp^{n-1}_{-a^2}$ denotes the simply connected
complete hyperbolic space of curvature $-a^2$. In each case,
\[
 I=[A,B]\quad(A<B),\qquad I=[A,\infty),
 \qquad\text{or}\qquad I=(-\infty,B],
\]
where each endpoint $A$ or $B$ is finite. The
sectional curvature is $0$ in \textup{(i)} and $-a^2$ in
\textup{(ii)}--\textup{(iii)}. The boundary has one or two components,
all noncompact.
\end{theorem}

The models in this theorem have concrete geometric descriptions.
The flat models are Euclidean half-spaces and slabs.
The exponential models are hyperbolic regions bounded by one or
two horospheres with a common ideal point; when there is only one
boundary component, the two half-line orientations correspond
to the two sides of a horosphere.
The hyperbolic-fiber models are regions bounded by one or two
leaves of the equidistant foliation about a totally geodesic
hyperplane.

The preceding results exhaust the possibilities in dimensions
$n\ge3$ and give the following global conclusions.

\begin{corollary}\label{cor:complete-global-topology}
Let $(M^n,g)$, $n\ge3$, be a complete connected
Riemannian manifold with nonempty smooth boundary. If
$\dim\ker\Phi^*=n$, then $M$ is simply connected and isometric
to one of the model domains listed in Theorems~\ref{thm:main},
\ref{thm:noncompact-main}, and~\ref{thm:noncompact-boundary-main}.
The boundary has one or two components, either all compact or all
noncompact. If $M$ is noncompact, it has exactly one end.
\end{corollary}

Conversely, every model in Theorems~\ref{thm:main},
\ref{thm:noncompact-main}, and~\ref{thm:noncompact-boundary-main}
satisfies $\dim\ker\Phi^*=n$.
For the models with connected boundary, this follows from
the kernel computations in \cite{Sheng2025}.
When the boundary has two components, one can choose $n$
linearly independent static potentials on $\Mspace_\kappa^n$
whose restrictions satisfy the boundary equation on both
components simultaneously, in accordance with the intersection
principle in \cite[Section~4.5]{Sheng2025}.
We include brief verifications in the proofs for completeness.

The remaining case is dimension two.
Although a single static potential already forces constant
Gaussian curvature, maximal boundary kernel dimension still
allows nontrivial global identifications.
Theorem~\ref{thm:surface-classification} completes the classification
by describing all complete surfaces with $\dim\ker\Phi^*=2$,
without assuming orientability.
In addition to the corresponding model domains, the classification
includes finite and universal coverings of concentric annuli
and geodesic-ball exteriors, as well as spherical M\"obius bands.

The maximality assumption in these classification results is
sharp for constant-curvature rigidity when $n\ge3$ and for
the bound on the number of boundary components when $n\ge2$.
Indeed, for every $n\ge3$,
Example~\ref{ex:submaximal-product} gives a compact product
with connected totally geodesic boundary and
$\dim\ker\Phi^*=n-1$ whose sectional curvature is not constant.
Moreover, for every $n\ge2$,
Example~\ref{ex:submaximal-many-components} gives compact
Euclidean domains with $\dim\ker\Phi^*=n-1$ and arbitrarily
many boundary components.

We now outline the global construction underlying these results.
In dimensions $n\ge3$, completeness allows the pair fields
\eqref{eq:intro-Xuv} to generate isometric actions preserving
the boundary.
A compact boundary component determines an effective
$\mathrm{SO}(n)$-action, while the case in which all boundary
components are noncompact leads to Euclidean or hyperbolic
symmetries.
The decisive step is to use the equivariant potential map
to determine the full stabilizer groups.
This identifies the regular orbits with the model spheres,
Euclidean spaces, or hyperbolic spaces themselves, excluding
nontrivial quotients of these orbits.

A global transverse section then removes possible twisting
and gives warped product coordinates away from possible centers.
The curvature identity determines the warping functions,
while the boundary equation and completeness determine the
admissible boundary leaves and endpoints.
Together, these steps resolve the global arrangement of the
orbits and complete the higher-dimensional classification.
For surfaces, the remaining task is to determine the angular
periods and reflections that preserve both potentials.

The paper is organized as follows.
Section~\ref{sec:static-equations} establishes the local curvature
theorem and the algebraic identities for static potentials
and their pair fields.
Section~\ref{sec:maximal-boundary} treats the boundary geometry,
complete flows, and the common global product construction.
Section~\ref{sec:orbit-space} determines the spherical orbits
and their endpoints.
Section~\ref{sec:compact-boundary-classification} proves the compact classification
and then the exterior classification.
Section~\ref{sec:noncompact-boundary-classification} treats the
models with noncompact boundary, and
Section~\ref{sec:surface-classification} completes the classification
for surfaces.
Finally, Section~\ref{sec:sharpness} gives the sharpness examples.

\section*{Acknowledgements}
We thank Richard Schoen for suggesting this problem and Feng Wang
for helpful discussions. Research of Yi was supported by NSFC (No.12501071) and NSFC (No.12531002).


\section{Local rigidity and the algebra of static potentials}
\label{sec:static-equations}

We first fix the notation and conventions used throughout the paper,
then establish the local identities underlying the classification.
The arguments in this section require no completeness assumption.

Let $(M^n,g)$, $n\ge2$, be a connected smooth Riemannian
manifold with nonempty smooth boundary.
Throughout the paper, completeness is assumed only when
explicitly stated. For a manifold with boundary, completeness
refers to the intrinsic Riemannian distance, including the
boundary.
We denote the entire boundary of $M$ by $\partial M$ and a
selected connected boundary component by $\Sigma$.
The notation $\Omega\subset\Mspace_\kappa^n$ is reserved for
concrete domains in model spaces, where $\Mspace_\kappa^n$ denotes the simply connected complete
space form of dimension $n$ and sectional curvature $\kappa$.

Let $D$ denote the Levi--Civita connection of $g$.
We also write $Du$ for the gradient of a function $u$ on $M$.
The divergence, Hessian, and Laplace--Beltrami operators of $g$
are denoted by $\operatorname{div}$, $\Hess$, and $\Delta$,
respectively, with
\[
 \Hess u=Ddu,
 \qquad
 \Delta u=\operatorname{tr}_g\Hess u.
\]
Our convention for the curvature tensor is
\begin{equation*}
 \operatorname{Rm}(X,Y)Z
 =D_XD_YZ-D_YD_XZ-D_{[X,Y]}Z.
\end{equation*}
Thus, for orthonormal vectors $X,Y$, the sectional curvature is
\[
 K(X\wedge Y)=g(\operatorname{Rm}(X,Y)Y,X).
\]
We write $R$ and $\Ric$ for the scalar curvature and
Ricci tensor of $g$, respectively.

Let $\nu$ denote the outward unit normal and $\hat g$ the
induced metric on $\partial M$.
Our convention for the second fundamental form $h$ and
mean curvature $H$ is
$$
h(X,Y)=-g(\nu,D_XY)=g(D_X\nu,Y),
 \quad H=\operatorname{tr}_{\hat g}h,
$$
where $X,Y\in T(\partial M)$.
The Levi--Civita connection of $\hat g$ is denoted by $\nabla$.
For a function $\phi$ on $\Sigma$, we write
$\nabla^\Sigma\phi$ for its gradient, $\Hess_{\hat g}\phi$
for its Hessian, and
$\Delta_\Sigma\phi=\Delta_{\hat g}\phi$ for its Laplacian
with respect to the induced metric.
Finally, for $u\in C^\infty(M)$, we write
$u_\nu=du(\nu)=D_\nu u$ for its derivative in the direction
of the outward unit normal.

\bigskip

We now turn to the static equations.
On a b-static manifold, the boundary is totally umbilic,
so the boundary equation \eqref{static bdry} is equivalent
to $u_\nu=\frac{H}{n-1}u$.
Taking the trace of \eqref{static int} and substituting
the resulting identity back into the equation gives an
equivalent Hessian formulation.
Consequently, a function $u\in C^\infty(M)$ belongs to
$\ker\Phi^*$ if and only if
\begin{equation}\label{eq:Hessian-static}
 \left\{
 \begin{aligned}
 \Hess u&=
 u\left(\Ric-\tfrac{R}{n-1}g\right)
 &&\text{in }M,\\
 u_\nu&=\tfrac{H}{n-1}u
 &&\text{on }\partial M.
 \end{aligned}
 \right.
\end{equation}

The interior Hessian equation shows that a static potential
is determined by its value and differential at one point.
For a b-static potential, the boundary equation also determines
the normal derivative from the boundary value.
Accordingly, for $p\in M$ and $q\in\Sigma$, define the
one-jet maps
\begin{align*}
 J_p:\ker L^*&\longrightarrow\R\oplus T_p^*M,
 &u&\longmapsto\bigl(u(p),du_p\bigr),\\
 J_q^\Sigma:\ker\Phi^*&\longrightarrow\R\oplus T_q^*\Sigma,
 &u&\longmapsto\bigl(u(q),d(u|_\Sigma)_q\bigr).
\end{align*}
The injectivity of these maps follows from the local
uniqueness argument in \cite[Corollary~2.4]{Corvino2000} and
\cite[Lemma~2.17 and Proposition~3.4]{ShengLocalized2024}.
The injectivity of restriction to a connected boundary
component is also established in
\cite[Proposition~2.5]{ShengZhao2026}, using the same
uniqueness principle.
These arguments use only local geodesics and connectedness,
and hence require no completeness assumption.
The dimension bounds and the isomorphism assertion below
then follow immediately.
We summarize these facts in the following proposition.

\begin{proposition}\label{prop:one-jet}
The map $J_p$ is injective for every $p\in M$.
For every connected boundary component $\Sigma$,
the map $J_q^\Sigma$ is injective for every $q\in\Sigma$,
and restriction $u\mapsto u|_\Sigma$ is injective on
$\ker\Phi^*$.
In particular,
\begin{equation*}
 \dim\ker L^*\le n+1,
 \qquad
 \dim\ker\Phi^*\le n.
\end{equation*}
If $\dim\ker\Phi^*=n$, then $J_q^\Sigma$ is an isomorphism
at every point of the boundary.
\end{proposition}

We next record the differential identity for the pair fields
\eqref{eq:intro-Xuv}; see
\cite[Proposition~3.2]{Lafontaine2009} and
\cite[Lemma~3.1]{HePetersenWylieRigidity2015}.

\begin{proposition}\label{prop:pair-Killing}
If $u,v\in\ker L^*$, then $X_{u,v}$ is Killing.
If $u,v\in\ker\Phi^*$, then $X_{u,v}$ is also tangent to
$\partial M$.
\end{proposition}

\begin{proof}
For $u,v\in\ker L^*$, the interior equation in
\eqref{eq:Hessian-static} gives
\[
 \begin{aligned}
 D_ZX_{u,v}
 &=du(Z)Dv-dv(Z)Du+uD_ZDv-vD_ZDu\\
 &=du(Z)Dv-dv(Z)Du
 \end{aligned}
\]
for every vector field $Z$.
Thus, the endomorphism $Z\mapsto D_ZX_{u,v}$ is skew-adjoint.
By the formula for the Lie derivative of the metric,
\[
 (\Lie_{X_{u,v}}g)(Y,Z)
 =g(D_YX_{u,v},Z)+g(D_ZX_{u,v},Y)=0.
\]
Therefore, $X_{u,v}$ is Killing.

If $u,v\in\ker\Phi^*$, the boundary equation in
\eqref{eq:Hessian-static} gives
\[
 \inner{X_{u,v}}{\nu}
 =uv_\nu-vu_\nu
 =\tfrac{H}{n-1}(uv-vu)=0
 \quad\text{on }\partial M.
\]
Thus, $X_{u,v}$ is tangent to $\partial M$.
\end{proof}

The preceding propositions give the dimension bounds for
static potentials and describe the Killing fields they generate.
The following theorem shows that $n$ independent static potentials force constant sectional curvature, although the one-jet bound is $n+1$.

\begin{theorem}\label{thm:interior-maximality}
\label{thm:interior-dimension-gap}
Let $(M^n,g)$, $n\ge3$, be a connected Riemannian manifold,
possibly with smooth boundary.
If $\dim\ker L^*\ge n$, then $g$ has constant sectional curvature $\kappa=R/[n(n-1)]$.
\end{theorem}

\begin{proof}
Since $\ker L^*\ne\{0\}$, the scalar curvature $R$ is constant \cite{FischerMarsden1975}.
Choose an $n$-dimensional subspace $W\subset\ker L^*$ and
$0\ne u_*\in W$.
By Proposition~\ref{prop:one-jet}, $u_*$ cannot vanish on a
nonempty open set. Hence $\{u_*\ne0\}\cap\operatorname{int}M$
is dense in $M$.
So it suffices to prove the curvature identity at a point $p$ of this set.

Set
\[
 W_p:=\{u\in W:u(p)=0\}.
\]
Since evaluation at $p$ is a nonzero linear functional on $W$,
the rank--nullity theorem gives $\dim W_p=n-1$.
One-jet uniqueness makes the gradient map on $W_p$ injective.
Its image therefore satisfies
\[
 E:=\{Du(p):u\in W_p\}\subset T_pM,\qquad \dim E=n-1.
\]

Set
\[
 B:=\Ric-\tfrac{R}{n-1}g,
 \qquad g(B^\sharp X,Y)=B(X,Y).
\]
We first show $B=-\kappa g$ at $p$.
The static equation gives $D_XDu=uB^\sharp X$, so
\[
 D_XD_YDu
 =du(X)B^\sharp Y+u(D_XB^\sharp)Y+uB^\sharp(D_XY).
\]
Interchanging $X,Y$ and using $D_XY-D_YX=[X,Y]$, we obtain
\[
 \operatorname{Rm}(X,Y)Du
 =du(X)B^\sharp Y-du(Y)B^\sharp X+u\bigl((D_XB^\sharp)Y-(D_YB^\sharp)X\bigr).
\]
For $u\in W_p$, we have $u(p)=0$, so the last term vanishes at $p$.
Thus
\begin{equation}\label{eq:interior-curvature-on-gradient-space}
 \operatorname{Rm}(X,Y)V
 =g(V,X)B^\sharp Y-g(V,Y)B^\sharp X
 \qquad(V\in E).
\end{equation}
Choose a unit normal $e_0$ to $E$ and an orthonormal basis
$e_1,\ldots,e_{n-1}$ diagonalizing $B|_{E\times E}$.
Taking $X=e_0$ and $Y=V=e_a$ in
\eqref{eq:interior-curvature-on-gradient-space} gives
\[
 \operatorname{Rm}(e_0,e_a)e_a=-B^\sharp e_0,
 \qquad a=1,\ldots,n-1.
\]
Since $\operatorname{Rm}(e_0,e_a)$ is skew-adjoint,
\[
 B(e_0,e_a)=-g(\operatorname{Rm}(e_0,e_a)e_a,e_a)=0.
\]
Thus $B$ is diagonal in this basis.
Write $B(e_0,e_0)=c$ and $B(e_a,e_a)=b_a$.
For distinct $a,b\in\{1,\ldots,n-1\}$,
\eqref{eq:interior-curvature-on-gradient-space} gives
\[
 -b_a=K(e_a\wedge e_b)=K(e_b\wedge e_a)=-b_b.
\]
Since $n-1\ge2$, all $b_a$ have a common value $b$.
The same identity gives $K(e_0\wedge e_a)=-c$.
Using the definition of $B$, we obtain
\[
 \begin{aligned}
 c+\tfrac{R}{n-1}&=\Ric(e_0,e_0)=\sum\nolimits_{a=1}^{n-1}K(e_0\wedge e_a)=-(n-1)c,\\
 b+\tfrac{R}{n-1}&=\Ric(e_a,e_a)=\sum\nolimits_{b\ne a} K(e_a\wedge e_b) + K(e_a\wedge e_0)=-(n-2)b-c.
 \end{aligned}
\]
Solving first for $c$ and then for $b$ yields
\[
 c=-\tfrac{R}{n(n-1)}=-\kappa,
 \qquad (n-1)b=-c-\tfrac{R}{n-1}=-(n-1)\kappa.
\]
This proves $B=-\kappa g$ at $p$.

To recover the full curvature tensor, define the skew-adjoint endomorphism
\[
 S_{X,Y}(Z):=\operatorname{Rm}(X,Y)Z
 -\kappa\bigl(g(Y,Z)X-g(X,Z)Y\bigr).
\]
By \eqref{eq:interior-curvature-on-gradient-space}, it vanishes
on $E$. Therefore, for every $V\in E$,
\[
 g(S_{X,Y}e_0,V)=-g(e_0,S_{X,Y}V)=0,
 \qquad g(S_{X,Y}e_0,e_0)=0.
\]
It follows that $S_{X,Y}e_0=0$, and hence $S_{X,Y}=0$.
We have proved
\begin{equation*}
 \operatorname{Rm}(X,Y)Z
 =\kappa\bigl(g(Y,Z)X-g(X,Z)Y\bigr)
\end{equation*}
at $p$. By density and continuity, the same identity holds on
all of $M$.
\end{proof}

\begin{remark}\label{rem:interior-dimension-two}
For $n=2$, one static potential already forces constant
Gaussian curvature, since $R$ is constant and $K=R/2$.
\end{remark}

When $g$ has constant curvature $\kappa$,
the system \eqref{eq:Hessian-static} is equivalent to
\begin{equation}\label{eq:space-form-static-system}
 \left\{
 \begin{aligned}
 \Hess u&=-\kappa ug
 &&\text{in }M,\\
 u_\nu&=\tfrac{H}{n-1}u
 &&\text{on }\partial M.
 \end{aligned}
 \right.
\end{equation}
In particular, this reduction applies under the assumption $\dim\ker\Phi^*=n$.

The reduced interior equation yields a bilinear form that
determines both the action of the pair fields on the potentials
and their Lie brackets; see
\cite[Section~8]{HePetersenWylieRigidity2015}.
We will use these identities in each of the global cases.

\begin{lemma}\label{lem:potential-bilinear-form}
Let $(M^n,g)$, $n\ge2$, be a connected Riemannian manifold
of constant sectional curvature $\kappa$, possibly with smooth
boundary.
For $u,v\in\ker L^*$, the quantity
\begin{equation}\label{eq:potential-bilinear-form}
 Q(u,v)=\inner{Du}{Dv}+\kappa uv
\end{equation}
is constant on $M$ and defines a symmetric bilinear form
on $\ker L^*$.
Moreover, for $u,v,w,z\in\ker L^*$, we have
\begingroup
\renewcommand{\mintagsep}{2pt}
\medmuskip=2mu
\thickmuskip=3mu
\begin{align}
 X_{u,v}(w)
 &= Q(v,w)u-Q(u,w)v,
 \label{eq:pair-action-bilinear-form}\\
 \lbrack X_{u,v},X_{w,z}\rbrack
 &= Q(v,w)X_{u,z}-Q(u,w)X_{v,z}-Q(v,z)X_{u,w}+Q(u,z)X_{v,w}.
 \label{eq:pair-bracket-bilinear-form}
\end{align}
\endgroup
\end{lemma}

\begin{proof}
Using the interior equation in
\eqref{eq:space-form-static-system}, we have, for every vector field $Z$,
\[
 Z(Q(u,v))
 =-\kappa u\,dv(Z)-\kappa v\,du(Z)
   +\kappa\bigl(v\,du(Z)+u\,dv(Z)\bigr)=0.
\]
Since $M$ is connected, $Q(u,v)$ is constant on $M$.
Substituting \eqref{eq:potential-bilinear-form} into
\[
 X_{u,v}(w)
 =u\inner{Dv}{Dw}-v\inner{Du}{Dw}
\]
gives \eqref{eq:pair-action-bilinear-form}.

To compute the Lie bracket, recall that
$[Y,Dw]=D(Yw)$ for every Killing field $Y$.
Applying this identity with $Y=X_{u,v}$, we obtain
\[
 [X_{u,v},X_{w,z}]
 =X_{X_{u,v}(w),z}+X_{w,X_{u,v}(z)}.
\]
Using \eqref{eq:pair-action-bilinear-form} and the constancy
of $Q$, we conclude that
\eqref{eq:pair-bracket-bilinear-form} holds.
\end{proof}

\section{Boundary geometry and global group actions}
\label{sec:maximal-boundary}
We begin by determining the geometry of each boundary component
and the restrictions of the b-static potentials.
For a boundary component $\Sigma$, set
\[
 \lambda:=\left.\tfrac{H}{n-1}\right|_\Sigma,
 \qquad
 \mathcal V:=\{u|_\Sigma:u\in\ker\Phi^*\}.
\]
Recall that the existence of a b-static potential implies \cite{ShengLocalized2024,Sheng2025}
\begin{equation}\label{eq:boundary-umbilic}
 h=\lambda\hat g,\qquad d^\Sigma\lambda=0,\qquad
 \Ric(X,\nu)=0\quad\text{for every }X\in T\Sigma.
\end{equation}

Under the maximality assumption $\dim\ker\Phi^*=n$,
Proposition~\ref{prop:one-jet} shows that $\dim\mathcal V=n$
and that, at every $q\in\Sigma$, the map
\[
 \mathcal V\longrightarrow\R\oplus T_q^*\Sigma,
 \qquad
 \phi\longmapsto\bigl(\phi(q),d\phi_q\bigr),
\]
is an isomorphism. Together with interior rigidity, these facts determine the boundary
geometry and the full space of boundary restrictions.

\begin{lemma}\label{lem:boundary-Hessian-system}
Let $(M^n,g)$, $n\ge3$, be a connected Riemannian manifold
with nonempty smooth boundary and $\dim\ker\Phi^*=n$.
Then the induced metric on each boundary component $\Sigma$
has constant sectional curvature
\begin{equation}\label{eq:boundary-curvature-general}
 k_\Sigma=\kappa+\lambda^2.
\end{equation}
Every $\phi\in\mathcal V$ satisfies
\begin{equation}\label{eq:boundary-Hessian-equation}
 \Hess_{\hat g}\phi=-k_\Sigma\phi\hat g,
\end{equation}
and $\mathcal V$ is the full solution space:
\begin{equation*}
 \mathcal V
 =\{f\in C^\infty(\Sigma):
       \Hess_{\hat g}f=-k_\Sigma f\hat g\}.
\end{equation*}
\end{lemma}

\begin{proof}
Since $\ker\Phi^*\subseteq\ker L^*$,
Theorem~\ref{thm:interior-maximality} shows that $g$ has
constant sectional curvature $\kappa=R/[n(n-1)]$.
By \eqref{eq:boundary-umbilic}, we have $h=\lambda\hat g$
with $\lambda$ constant on $\Sigma$.
The Gauss equation therefore gives
\eqref{eq:boundary-curvature-general}.

Let $\phi=u|_\Sigma$ for $u\in\ker\Phi^*$.
The Gauss formula gives
\begin{equation}\label{eq:boundary-Gauss-Hessian}
 \Hess_{\hat g}\phi=(\Hess u)^T-u_\nu h.
\end{equation}
By \eqref{eq:space-form-static-system}, we have
$(\Hess u)^T=-\kappa\phi\hat g$ and $u_\nu=\lambda\phi$.
Substituting these identities and $h=\lambda\hat g$ into
\eqref{eq:boundary-Gauss-Hessian}, we obtain
\[
 \Hess_{\hat g}\phi
 =-\kappa\phi\hat g-\lambda^2\phi\hat g
 =-k_\Sigma\phi\hat g.
\]

The one-jet uniqueness argument underlying
Proposition~\ref{prop:one-jet} also applies to
\eqref{eq:boundary-Hessian-equation}, so its solution space has
dimension at most $1+\dim\Sigma=n$.
Since it contains the $n$-dimensional space $\mathcal V$,
the two spaces coincide.
\end{proof}

On a complete boundary component, the full Hessian
solution space excludes nontrivial space-form quotients.

\begin{lemma}\label{lem:maximal-Hessian-system}
\label{lem:complete-boundary-model}
Under the hypotheses of Lemma~\ref{lem:boundary-Hessian-system},
let $\Sigma$ be a boundary component whose induced metric
is complete. Then
\[
 (\Sigma,\hat g)\cong\Mspace_{k_\Sigma}^{n-1}.
\]
Under this isometry, $\mathcal V$ consists of all affine
functions on $\R^{n-1}$ when $k_\Sigma=0$, and of the
restrictions of ambient linear functions to the spherical
or hyperboloid model when $k_\Sigma>0$ or $k_\Sigma<0$,
respectively.

If $\Sigma$ is compact, then $k_\Sigma>0$,
\begin{equation}\label{eq:Hessian-system-rigidity}
 (\Sigma,\hat g)\cong\Sph^{n-1}(k_\Sigma^{-1/2}),
\end{equation}
and
\begin{equation}\label{eq:Hessian-system-eigenspace}
 \mathcal V
 =\ker\bigl(\Delta_{\hat g}+(n-1)k_\Sigma\bigr).
\end{equation}
\end{lemma}

\begin{proof}
Set $m=n-1$ and $k=k_\Sigma$.
By Lemma~\ref{lem:boundary-Hessian-system}, $\hat g$ has
constant sectional curvature $k$.
Completeness and the Killing--Hopf theorem allow us to
identify the universal Riemannian covering with
\[
 \pi:(\Mspace_k^m,\widetilde g)\longrightarrow(\Sigma,\hat g),
 \qquad \widetilde g=\pi^*\hat g.
\]
It remains to show that $\pi$ is injective.

Since $\pi$ is a surjective local isometry, $\pi^*\mathcal V$ is an
$(m+1)$-dimensional subspace of the solutions of
$\Hess_{\widetilde g}F=-kF\widetilde g$.
On the spherical and hyperboloid models, this solution space is
spanned by the ambient coordinate functions. On the Euclidean model,
it is $\Span\{1,x_1,\ldots,x_m\}$ \cite{Sheng2025}.
It therefore has dimension $m+1$, so
\begin{equation}\label{eq:Hessian-lifted-space}
 \pi^*\mathcal V
 =\{F\in C^\infty(\Mspace_k^m):
       \Hess_{\widetilde g}F=-kF\widetilde g\}.
\end{equation}

If $\pi(p)=\pi(q)$, every function in $\pi^*\mathcal V$
takes the same value at $p$ and $q$.
In particular, the functions listed above agree at these
two points. Since they separate points on each model,
we obtain $p=q$.
Thus, $\pi$ is injective and, being a bijective local isometry,
is a global isometry.
This proves the asserted identification of $(\Sigma,\hat g)$
and the description of $\mathcal V$.

If $\Sigma$ is compact, only the spherical model is possible.
Hence $k_\Sigma>0$ and
\eqref{eq:Hessian-system-rigidity} holds.
On this sphere, the coordinate functions span the full
eigenspace of $-\Delta_{\hat g}$ with eigenvalue
$(n-1)k_\Sigma$, which gives
\eqref{eq:Hessian-system-eigenspace}.
\end{proof}

These results together complete the proof of
Theorem~\ref{thm:local-main}.

To compare different boundary components, we examine
the restriction of $Q$ to $\ker\Phi^*$.
The next lemma shows its inertia determines
the sign of $k_\Sigma$.

\begin{lemma}\label{lem:potential-form-boundary}
Under the hypotheses of Lemma~\ref{lem:boundary-Hessian-system},
let $\Sigma$ be a boundary component. For all $u,v\in\ker\Phi^*$,
we have
\begin{equation}\label{eq:potential-form-boundary}
 Q(u,v)
 =\inner{d^\Sigma u}{d^\Sigma v}_{\hat g}+k_\Sigma uv
 \qquad\text{on }\Sigma.
\end{equation}
The inertia of $Q|_{\ker\Phi^*}$, in the order positive, negative,
and zero, is
\begin{equation}\label{eq:potential-form-inertia}
 \begin{cases}
  (n,0,0),   & k_\Sigma>0,\\
  (n-1,0,1), & k_\Sigma=0,\\
  (n-1,1,0), & k_\Sigma<0.
 \end{cases}
\end{equation}
In particular, the sign of $k_\Sigma$ is the same for all
boundary components.
\end{lemma}

\begin{proof}
Decomposing $Du$ and $Dv$ into their tangential and normal
components along $\Sigma$, we obtain
\[
 Q(u,v)
 =\inner{d^\Sigma u}{d^\Sigma v}_{\hat g}
  +u_\nu v_\nu+\kappa uv.
\]
Since $u_\nu=\lambda u$, $v_\nu=\lambda v$, and
$k_\Sigma=\kappa+\lambda^2$, this gives
\eqref{eq:potential-form-boundary}.

Fix $q\in\Sigma$. By Proposition~\ref{prop:one-jet}, the boundary
one-jet map $J_q^\Sigma$ identifies $Q|_{\ker\Phi^*}$ with the
bilinear form on $\R\oplus T_q^*\Sigma$ given by
\[
 ((s,\xi),(t,\eta))
 \longmapsto k_\Sigma st+\inner{\xi}{\eta}_{\hat g_q}.
\]
Choosing an orthonormal basis of $T_q^*\Sigma$ gives the matrix
$\operatorname{diag}(k_\Sigma,1,\ldots,1)$, which yields
\eqref{eq:potential-form-inertia}.
Since $Q|_{\ker\Phi^*}$ is independent of $\Sigma$, the sign of
$k_\Sigma$ is also independent of $\Sigma$.
\end{proof}

If $M$ is complete, each boundary component is complete
with respect to its induced metric.
Combining Lemmas~\ref{lem:complete-boundary-model} and \ref{lem:potential-form-boundary}, we find that all boundary components
are simply connected space forms whose sectional curvatures
have \emph{the same sign}.
Together with \eqref{eq:boundary-curvature-general},
this gives the following corollary.

\begin{corollary}\label{cor:complete-noncompact-boundary}
Under the hypotheses of Lemma~\ref{lem:boundary-Hessian-system},
assume that $M$ is complete.
Then the boundary components are either all round spheres,
all isometric to $\R^{n-1}$, or all isometric to simply
connected hyperbolic spaces.
If some boundary component is noncompact, then $\kappa\le0$,
with strict inequality in the hyperbolic case.
\end{corollary}

We next pass from the potential algebra to global isometries.
The following lemma combines completeness of boundary-tangent
Killing fields with the Lie--Palais theorem.

\begin{lemma}\label{lem:complete-Killing-integration}
Let $(M,g)$ be a complete manifold with smooth boundary.
\begin{enumerate}[label=\textup{(\roman*)}]
\item Every smooth Killing field tangent to $\partial M$ is complete.
Its flow preserves both $\partial M$ and $\operatorname{int}M$.

\item Let $G$ be a simply connected Lie group with Lie algebra
$\mathfrak g$, and let $\iota$ be a Lie algebra homomorphism from
$\mathfrak g$ into the Killing fields tangent to $\partial M$.
Denote the flow of $\iota(A)$ by $\varphi_t^A$.
Then there is a unique smooth isometric left $G$-action on $M$
satisfying
\[
 \exp(-tA)\cdot p=\varphi_t^A(p)
 \qquad
 (A\in\mathfrak g,\ t\in\R,\ p\in M).
\]
This action preserves $\partial M$.
\end{enumerate}
\end{lemma}

\begin{proof}
Let $Y$ be a smooth Killing field tangent to $\partial M$.
By the flow theorem for manifolds with boundary
\cite[Theorem~9.34]{Lee2013}, its local flow preserves both
$\partial M$ and $\operatorname{int}M$.
Along any maximal integral curve $\gamma$, the Killing equation gives
\[
 \tfrac{d}{dt}|Y|^2(\gamma(t))=2g(D_YY,Y)=0.
\]
Thus $\gamma$ has constant speed $c$, and
\[
 d_M(\gamma(s),\gamma(t))\le c|s-t|.
\]
At a finite endpoint of its maximal interval, $\gamma$ would
converge to a point of $M$ by completeness and extend by the
local flow theorem. Hence $Y$ is complete, proving \textup{(i)}.

By \textup{(i)}, every field $\iota(A)$ is complete and its flow
preserves $\operatorname{int}M$.
The Lie--Palais theorem \cite{Palais1957} gives a smooth left
$G$-action $\alpha$ on $\operatorname{int}M$ satisfying
\begin{equation}\label{eq:complete-integration-interior}
 \alpha(\exp(-tA),p)=\varphi_t^A(p).
\end{equation}
Since $G$ is connected, every $a\in G$ has a finite factorization
\[
 a=\exp(-t_1A_1)\cdots\exp(-t_rA_r).
\]
Define
\begin{equation}\label{eq:complete-integration-extension}
 F_a=\varphi_{t_1}^{A_1}\circ\cdots\circ\varphi_{t_r}^{A_r}
 :M\to M.
\end{equation}
On $\operatorname{int}M$, $F_a=\alpha(a,\cdot)$ by the action
property and \eqref{eq:complete-integration-interior}.
Since $\operatorname{int}M$ is dense in $M$, continuity shows
that $F_a$ is independent of the factorization and that
\[
 F_e=\operatorname{id}_M,\qquad F_{ab}=F_a\circ F_b.
\]

Choose a basis $B_1,\ldots,B_d$ of $\mathfrak g$.
In coordinates of the second kind near $a$, we have
\[
 \begin{aligned}
 a'&=a\exp(-s_1B_1)\cdots\exp(-s_dB_d),\\
 F_{a'}&=F_a\circ\varphi_{s_1}^{B_1}\circ\cdots
               \circ\varphi_{s_d}^{B_d}.
 \end{aligned}
\]
The second identity proves that $(a,p)\mapsto F_a(p)$ is
smooth up to $\partial M$.
Taking a single exponential in
\eqref{eq:complete-integration-extension} gives
$F_{\exp(-tA)}=\varphi_t^A$.
Each $F_a$ is a boundary-preserving isometry, being a
composition of Killing flows.
Uniqueness follows from uniqueness of flows, since
exponentials generate $G$.
For the boundary-tangent formulation, see also
\cite[Section~1.2, pp.~220--221]{Hirsch2011}.
\end{proof}

We also record the product construction used below.
It gives a global product from an equivariant map whose
point stabilizers agree with those of the target.
Here $G$ acts on $G/K$ by left translation, and $G_x$ denotes
the stabilizer of a point $x$ in either $M$ or $G/K$.

\begin{lemma}\label{lem:equivariant-product}
Let a Lie group $G$ act smoothly on a nonempty manifold $M$,
possibly with smooth boundary.
Let $K\subset G$ be a closed subgroup and
$\Theta:M\to G/K$ a smooth $G$-equivariant map satisfying
\[
 G_p=G_{\Theta(p)}
 \qquad\text{for every }p\in M.
\]
Set $C=\Theta^{-1}(eK)$, where $e$ is the identity of $G$.
Then $C$ is a closed embedded submanifold transverse to
$\partial M$, with $\partial C=C\cap\partial M$, and $K$
fixes $C$ pointwise.
The map
\begin{equation}\label{eq:equivariant-product}
 F:C\times G/K\longrightarrow M,
 \qquad F(c,gK)=g\cdot c,
\end{equation}
is a $G$-equivariant diffeomorphism of manifolds with boundary
satisfying
\[
 \Theta(F(c,gK))=gK,
\]
where $G$ acts on the product by
$a\cdot(c,gK)=(c,agK)$.
\end{lemma}

\begin{proof}
Equivariance and transitivity imply that $\Theta$ is onto
and that its differential maps the orbit directions onto
$T_{\Theta(p)}(G/K)$ at every $p\in M$.
At boundary points these directions lie in $T_p\partial M$,
so both $\Theta$ and $\Theta|_{\partial M}$ are submersions.
Since $G/K$ is Hausdorff, $C$ is closed.
The regular-level theorem for manifolds with boundary shows
that $C$ is an embedded submanifold transverse to $\partial M$,
with
\[
 \dim C=\dim M-\dim(G/K),
 \qquad
 \partial C=C\cap\partial M.
\]

For $c\in C$, the stabilizer hypothesis gives $G_c=K$.
Thus $K$ fixes $C$ pointwise, and $F$ is well defined.
Choose a smooth local section $\sigma:U\to G$ of $G\to G/K$.
The restriction of $F$ to $C\times U$ and its inverse on
$\Theta^{-1}(U)$ are given by
\[
 \begin{aligned}
 F(c,y)&=\sigma(y)\cdot c,\\
 F^{-1}(p)&=
 \bigl(\sigma(\Theta(p))^{-1}\cdot p,\Theta(p)\bigr).
 \end{aligned}
\]
Equivariance ensures that the first component of the inverse
lies in $C$.
These mutually inverse local formulas are smooth up to the
boundary and show that $F$ is a global diffeomorphism.
Its equivariance and the identity
$\Theta(F(c,gK))=gK$ follow directly from its definition.
\end{proof}

Now let $(M^n,g)$, $n\ge3$, be complete and connected,
with $\dim\ker\Phi^*=n$ and a compact boundary component $\Sigma$.
By Lemma~\ref{lem:complete-boundary-model}, identify $\Sigma$ with
$\Sph^{n-1}(r_0)\subset\R^n$. Its ambient coordinate functions
$x_i=r_0\theta_i$ are the restrictions of a unique basis
$u_1,\ldots,u_n$ of $\ker\Phi^*$:
\begin{equation}\label{eq:normalized-boundary-basis}
 u_i|_\Sigma=x_i,\qquad i=1,\ldots,n.
\end{equation}
Set
\begin{equation}\label{eq:Xij-definition}
 X_{ij}=u_iDu_j-u_jDu_i,
 \qquad X_{ji}=-X_{ij},\qquad X_{ii}=0,
\end{equation}
and define the potential map
\begin{equation}\label{eq:potential-map}
 \mathbf u=(u_1,\ldots,u_n):M\longrightarrow\R^n.
\end{equation}

\begin{proposition}\label{prop:rotation-action}
Let $(M^n,g)$, $n\ge3$, be complete and connected, with nonempty
smooth boundary and $\dim\ker\Phi^*=n$. Let $\Sigma$ be a
compact boundary component and choose the basis
\eqref{eq:normalized-boundary-basis}. Then the fields $X_{ij}$
integrate to an effective isometric $\mathrm{SO}(n)$-action preserving
every boundary component. Moreover,
\begin{equation}\label{eq:rotation-on-potentials}
 X_{ij}(u_k)=\delta_{jk}u_i-\delta_{ik}u_j,
\end{equation}
and the potential map is equivariant:
\begin{equation}\label{eq:potential-equivariance}
 \mathbf u(A\cdot p)=A\mathbf u(p),
 \qquad A\in\mathrm{SO}(n),\quad p\in M.
\end{equation}
\end{proposition}

\begin{proof}
By \eqref{eq:potential-form-boundary} and the
coordinate identities on $\Sph^{n-1}(r_0)$, we have
\begin{equation}\label{eq:normalized-potential-form}
 Q(u_i,u_j)
 =\inner{d^\Sigma x_i}{d^\Sigma x_j}_{\hat g}
   +r_0^{-2}x_ix_j=\delta_{ij}.
\end{equation}
On $M$, Lemma~\ref{lem:potential-bilinear-form} gives
\eqref{eq:rotation-on-potentials} and
\begin{equation}\label{eq:so-brackets}
 [X_{ij},X_{kl}]
 =\delta_{jk}X_{il}-\delta_{ik}X_{jl}
  -\delta_{jl}X_{ik}+\delta_{il}X_{jk}.
\end{equation}

For a left action we use the fundamental-field convention
\begin{equation}\label{eq:rotation-fundamental-convention}
 A_M(p)=\left.\tfrac{d}{dt}\right|_{t=0}\exp(-tA)\cdot p,
\end{equation}
so $A\mapsto A_M$ is a Lie algebra homomorphism. The standard
matrices $E_{ij}=e_ie_j^{\mathsf T}-e_je_i^{\mathsf T}$ satisfy
\eqref{eq:so-brackets}, and hence $E_{ij}\mapsto X_{ij}$ defines a
homomorphism from $\mathfrak{so}(n)$ to boundary-tangent Killing
fields. Lemma~\ref{lem:complete-Killing-integration} integrates it to
an isometric action $\Psi$ of $\operatorname{Spin}(n)$.

Let $\varpi:\operatorname{Spin}(n)\to\mathrm{SO}(n)$ be the
standard covering. On $\Sigma$,
\begin{equation}\label{eq:rotation-boundary-fields}
 X_{ij}|_\Sigma
 =x_i\nabla^\Sigma x_j-x_j\nabla^\Sigma x_i,
\end{equation}
so uniqueness of integration gives
$\Psi(s)(q)=\varpi(s)q$ for $q\in\Sigma$.
Faithfulness of the standard sphere action implies
$\ker\Psi\subseteq\ker\varpi$. Conversely, if $s\in\ker\varpi$,
the isometry $\Psi(s)$ fixes $\Sigma$ pointwise and preserves its
inward unit normal. It fixes an inward collar by normal geodesics,
and hence is the identity on the connected manifold $M$.
Thus $\ker\Psi=\ker\varpi$, and the action descends to an effective
$\mathrm{SO}(n)$-action. Connectedness of this group implies that
each boundary component is preserved.

Finally, \eqref{eq:rotation-on-potentials} is equivalent to
\begin{equation}\label{eq:potential-infinitesimal-equivariance}
 d\mathbf u(A_M)=-A\mathbf u,
 \qquad A\in\mathfrak{so}(n).
\end{equation}
Integrating along $\exp(-tA)\cdot p$ gives
$\mathbf u(\exp(-tA)\cdot p)=\exp(-tA)\mathbf u(p)$.
Since exponentials generate $\mathrm{SO}(n)$, this proves
\eqref{eq:potential-equivariance}.
\end{proof}


\section{The equivariant potential map and the orbit space}
\label{sec:orbit-space}

Throughout this section, let $(M^n,g)$, $n\ge3$, be complete and
connected, with nonempty smooth boundary and $\dim\ker\Phi^*=n$.
Assume that one boundary component $\Sigma$ is compact. Choose the
normalized basis in \eqref{eq:normalized-boundary-basis}, and let
$G=\mathrm{SO}(n)$ act as in Proposition~\ref{prop:rotation-action}.
The same construction will treat compact manifolds and complete
noncompact manifolds with a compact boundary component.

Write $G_p$ for the stabilizer of $p\in M$ and $G_z$ for the
stabilizer of $z\in\R^n$ under the standard action. We first identify
the full stabilizers, then apply Lemma~\ref{lem:equivariant-product}
to obtain global coordinates away from the fixed points.

\subsection{Exact isotropy groups}

\begin{proposition}\label{prop:isotropy-rigidity}
Under the standing assumptions, the following hold:
\begin{enumerate}[label=\textup{(\roman*)}]
\item A point $p\in M$ is fixed by $G$ if and only if $\mathbf u(p)=0.$
\item At every nonfixed point,
\begin{equation}\label{eq:exact-isotropy}
 G_p=G_{\mathbf u(p)}\cong\mathrm{SO}(n-1),
 \qquad G\cdot p\cong\Sph^{n-1}.
\end{equation}
\item Every fixed point is isolated and interior, with isotropy
representation equivalent to the standard $\mathrm{SO}(n)$-action
on $\R^n$.
\end{enumerate}
\end{proposition}

\begin{proof}
If $\mathbf u(p)=0$, all the generating fields $X_{ij}$ vanish at
$p$ by \eqref{eq:Xij-definition}, so the connected group $G$ fixes $p$. Conversely, equivariance \eqref{eq:potential-equivariance}
implies that $\mathbf u(p)$ is fixed by every standard rotation
whenever $p$ is fixed. Thus $\mathbf u(p)=0$.

Suppose that $\mathbf u(p)\ne0$. By \eqref{eq:potential-equivariance}, 
$G_p\subseteq G_{\mathbf u(p)}$. After rotating $p$, we may assume
$\mathbf u(p)=ae_1$, where $a>0$. The stabilizer inclusion gives
\begin{equation}\label{eq:orbit-lower-bound}
 \dim(G\cdot p)\ge\dim G-\dim\mathrm{SO}(n-1)=n-1.
\end{equation}
On the other hand, $u_i(p)=0$ for $i\ge2$, so
\begin{equation}\label{eq:orbit-upper-bound}
 T_p(G\cdot p)=\Span\{X_{12}(p),\ldots,X_{1n}(p)\},
 \qquad \dim(G\cdot p)\le n-1.
\end{equation}
Hence $G_p$ and $G_{\mathbf u(p)}$ have the same dimension.
Since $G_p$ is closed and $G_{\mathbf u(p)}\cong\mathrm{SO}(n-1)$
is connected,
\begin{equation}\label{eq:isotropy-component-inclusions}
 G_{\mathbf u(p)}=(G_{\mathbf u(p)})^0=(G_p)^0
 \subseteq G_p\subseteq G_{\mathbf u(p)},
\end{equation}
where the superscript $0$ denotes the identity component. This proves
\textup{(ii)}.

At a fixed point, all potentials vanish. The boundary one-jet
isomorphism excludes such a point from $\partial M$.
One-jet uniqueness and $\dim\ker\Phi^*=n$ show that the map
$u\mapsto du_p$ is an isomorphism onto $T_p^*M$. Thus
$d\mathbf u_p$ is invertible. Differentiating \eqref{eq:potential-equivariance} gives
\begin{eqnarray*}
 d\mathbf u_p\circ d(A\cdot)_p=A\circ d\mathbf u_p
 \qquad\text{for every }A\in G,
\end{eqnarray*}
which identifies the isotropy representation with the standard one.
The inverse function theorem shows that $p$ is an isolated zero
of $\mathbf u$, proving \textup{(iii)}.
\end{proof}

\subsection{The orbit interval and global coordinates}

\begin{lemma}\label{lem:orbit-interval}
Under the standing assumptions, every boundary
component is a single principal orbit. With the quotient metric, the orbit space is
\begin{eqnarray*}
 M/G\cong
 \begin{cases}
  [0,L],&M\text{ compact},\\
  [0,\infty),&M\text{ noncompact},
 \end{cases}
\end{eqnarray*}
where $L>0$ and the selected component $\Sigma$ corresponds to $0$.
Its endpoints correspond precisely to boundary components or fixed
points. In particular, there are at most two boundary components.
In the noncompact case, $\partial M=\Sigma$ and there are no
fixed points.
\end{lemma}

\begin{proof}
 By Proposition~\ref{prop:isotropy-rigidity}, boundary points are
nonfixed, and all nonfixed orbits have stabilizers conjugate to
$\mathrm{SO}(n-1)$. Each boundary orbit is open in its component by
dimension and closed by compactness of $G$. Thus every boundary
component is a single orbit, and all nonfixed orbits are principal.

Let $\mathcal I=M/G$. Since $G$ is compact, $\mathcal I$ is
Hausdorff and the quotient map is proper
\cite[Chapter~I, Theorem~3.1]{Bredon1972}. At a nonfixed interior
point, the connected stabilizer acts orthogonally, hence trivially, on
the one-dimensional normal space. The equivariant tubular-neighborhood
theorem \cite[Chapter~VI, Theorem~2.2]{Bredon1972}, the inward normal
collar at a boundary orbit, and the standard isotropy representation at a fixed point give invariant neighborhoods $U$ with
\begin{eqnarray*}
 U/G\cong
 \begin{cases}
  (-\varepsilon,\varepsilon),
    &p\in\operatorname{int}M\text{ nonfixed},\\
  [0,\varepsilon),
    &p\in\partial M\text{ or }p\text{ fixed}.
 \end{cases}
\end{eqnarray*}
Hence $\mathcal I$ is a connected one-dimensional manifold whose
boundary consists precisely of the boundary and fixed-point orbits.

Equip $\mathcal I$ with
\begin{eqnarray*}
 d_{\mathcal I}(G\cdot p,G\cdot q)
 =\min_{A\in G}d_M(p,A\cdot q),
\end{eqnarray*}
where $d_M$ is the intrinsic length distance. The quotient
projection is $1$-Lipschitz, and projecting nearly minimizing curves
shows that $d_{\mathcal I}$ is a length metric. After shrinking the
local neighborhoods, the Gauss lemma identifies this metric with
$|s-t|$ in normal arclength coordinates.

To prove completeness, take a Cauchy sequence of orbits and pass to a
subsequence whose successive distances are less than $2^{-j}$. Choose
representatives $p_j$ recursively so that
\begin{eqnarray*}
 d_M(p_j,p_{j+1})
 =d_{\mathcal I}(G\cdot p_j,G\cdot p_{j+1})<2^{-j}.
\end{eqnarray*}
Then $p_j$ is Cauchy and converges in $M$. The subsequence of
orbits, and hence the original sequence, therefore converges.

Since $\mathcal I$ is a connected one-dimensional manifold with
nonempty boundary and $d_{\mathcal I}$ is a complete arclength
metric, it is isometric to a closed bounded interval or a closed
half-line. Properness of the quotient map shows that the first
case occurs exactly when $M$ is compact. Distinct
boundary components give distinct endpoints, and $\Sigma$ gives the
endpoint $0$. A closed interval has two endpoints, so there are at most
two boundary components; a half-line has only one, so in the noncompact
case $\partial M=\Sigma$ and there are no fixed points.
\end{proof}

Let $\pi:M\to\mathcal I$ be the quotient map. In the compact
case write $\mathcal I=[0,L]$; in the noncompact case set $L=\infty$
and write $\mathcal I=[0,L)$. Define
\begin{eqnarray*}
 M_* =\{p\in M:\mathbf u(p)\ne0\},
 \qquad M_{\mathrm{reg}}=\pi^{-1}((0,L)),
\end{eqnarray*}
and set
\begin{equation}\label{eq:principal-isotropy}
 K=G_{e_1}\cong\mathrm{SO}(n-1).
\end{equation}
We identify $G/K$ with the unit sphere by $AK\mapsto Ae_1$.

\begin{lemma}\label{lem:global-section}
Under the standing assumptions, the normalized potential map
\begin{eqnarray*}
 \Theta=\frac{\mathbf u}{|\mathbf u|}:M_*\longrightarrow\Sph^{n-1}
\end{eqnarray*}
gives a global product of $M_*$ with spherical orbits. More
precisely, there is a unit-speed geodesic $q:[0,L)\to M_*$ with
$\mathbf u(q(0))=r_0e_1$ such that
\begin{equation}\label{eq:horizontal-section}
 \pi(q(t))=t,\qquad
 \dot q(t)\perp T_{q(t)}(G\cdot q(t)).
\end{equation}
The map
\begin{equation}\label{eq:global-F}
 F:(0,L)\times G/K\longrightarrow M_{\mathrm{reg}},
 \qquad F(t,AK)=A\cdot q(t),
\end{equation}
is an equivariant diffeomorphism, and
\begin{equation}\label{eq:global-warped-and-potential}
 F^*g=dt^2+\rho(t)^2g_{\Sph^{n-1}},\qquad
 \mathbf u(F(t,\theta))=a(t)\theta
\end{equation}
for smooth positive functions $\rho,a$. The product extends smoothly
to every boundary endpoint, with $\rho(0)=a(0)=r_0$.
\end{lemma}

\begin{proof}
By \eqref{eq:potential-equivariance} and
Proposition~\ref{prop:isotropy-rigidity}, for every $A\in G$ and
$p\in M_*$,
\begin{eqnarray*}
 \Theta(A\cdot p)=A\Theta(p),\qquad G_p=G_{\Theta(p)}.
\end{eqnarray*}
Set $C=\Theta^{-1}(e_1)$. By
Lemma~\ref{lem:equivariant-product}, $K$ fixes $C$ pointwise, and
\begin{eqnarray*}
 C\times G/K\longrightarrow M_*,
 \qquad (c,AK)\longmapsto A\cdot c
\end{eqnarray*}
is an equivariant diffeomorphism, including the boundary.
Thus $C$ meets each nonfixed orbit exactly once, and
$\pi|_C$ identifies $C$ with the nonfixed part of $\mathcal I$.

Fix $c\in C$ and $v\in T_cC$, and let $v^\top$ be the orthogonal
projection of $v$ onto $T_c(G\cdot c)$. Since $K$ fixes $C$
pointwise and acts by isometries, for every $k\in K$,
\begin{eqnarray*}
 dk_c(v)=v,\qquad dk_c(v^\top)=(dk_c(v))^\top=v^\top.
\end{eqnarray*}
The restriction of $d\Theta_c$ to $T_c(G\cdot c)$ is a
$K$-equivariant isomorphism onto $e_1^\perp$. The standard
representation of $K=\mathrm{SO}(n-1)$ on $e_1^\perp$ has no
nonzero fixed vector, since $n\ge3$. Hence $v^\top=0$, and
\begin{eqnarray*}
 T_cC\perp T_c(G\cdot c).
\end{eqnarray*}

Parametrize $C$ by arclength from $\Sigma$. By the local quotient
metric in Lemma~\ref{lem:orbit-interval},
\begin{eqnarray*}
 |\dot q(t)|=1,\qquad
 (\pi\circ q)'(t)=1,\qquad
 \pi(q(0))=0.
\end{eqnarray*}
Thus $\pi(q(t))=t$, proving \eqref{eq:horizontal-section}.
The product diffeomorphism gives
\eqref{eq:global-F}.

For $\theta=Ae_1$, we have
$\partial_tF=dA_{q(t)}\dot q(t)$. Since $A$ is an isometry,
\begin{eqnarray*}
 (F^*g)(\partial_t,\partial_t)
 &=&|dA_{q(t)}\dot q(t)|^2=1,\\
 (F^*g)(\partial_t,V)&=&0
 \qquad\text{for every }V\in T_\theta\Sph^{n-1}.
\end{eqnarray*}
Each orbit metric is $G$-invariant. Since $K$ acts transitively
on the standard unit sphere in $e_1^\perp$, it is a positive
multiple of the round metric. Therefore
\begin{eqnarray*}
 F^*g=dt^2+\rho(t)^2g_{\Sph^{n-1}}
\end{eqnarray*}
for a smooth positive function $\rho$.
Since $g_{tt}=1$ and $g_{t\alpha}=0$, we have
$D_{\partial_t}\partial_t=0$. Thus $q$ is a unit-speed geodesic.

Set $a(t)=|\mathbf u(q(t))|>0$. Since $\Theta(q(t))=e_1$,
\begin{eqnarray*}
 \mathbf u(q(t))=a(t)e_1.
\end{eqnarray*}
For $\theta=Ae_1$, \eqref{eq:potential-equivariance} gives
\begin{eqnarray*}
 \mathbf u(F(t,\theta))
 =\mathbf u(A\cdot q(t))
 =A\mathbf u(q(t))
 =a(t)\theta.
\end{eqnarray*}
This proves \eqref{eq:global-warped-and-potential}.

The product diffeomorphism includes every boundary orbit,
and arclength is a smooth coordinate on $C$ up to its boundary.
Thus $F$, $\rho$, and $a$ extend smoothly to every boundary
endpoint. By \eqref{eq:normalized-boundary-basis},
\begin{eqnarray*}
 \mathbf u(q(0))=r_0e_1,\qquad
 \rho(0)=a(0)=r_0.
\end{eqnarray*}
\end{proof}

\subsection{Regularity at the endpoints}

At a finite endpoint $t_*$ of $\mathcal I$, set $s=|t-t_*|$ and
write $\rho_*(s)=\rho(t)$ near $t_*$.

\begin{lemma}\label{lem:endpoint-regularity}
Under the standing assumptions, the warping function has the following
behavior at each finite endpoint:
\begin{enumerate}[label=\textup{(\roman*)}]
\item At a boundary orbit, $\rho_*$ is smooth up to $s=0$ and
$\rho_*(0)>0$.
\item At a fixed point, $\rho_*$ is smooth up to $s=0$ and
\begin{eqnarray*}
 \rho_*(0)=0,\qquad \rho_*'(0)=1.
\end{eqnarray*}
\end{enumerate}
\end{lemma}

\begin{proof}
At a boundary orbit, \textup{(i)} follows from
Lemma~\ref{lem:global-section}: the product extends smoothly to the
boundary, where the orbit metric is nondegenerate.

Suppose that $t_*$ corresponds to a fixed point $p$.
Since $q$ is unit speed, $d_M(q(t),q(t'))\le|t-t'|$.
By completeness, $q(t)\to p_*$ as $t\to t_*$.
By \eqref{eq:horizontal-section},
\begin{eqnarray*}
 \pi(p_*)=\lim_{t\to t_*}\pi(q(t))=t_*,\qquad
 \pi^{-1}(t_*)=\{p\}.
\end{eqnarray*}
Thus $p_*=p$. Passing to the limit in the distance estimate and
using the $1$-Lipschitz property of $\pi$, we obtain
\begin{eqnarray*}
 s=|t-t_*|
 \le d_M(q(t),p)\le s.
\end{eqnarray*}
Hence the geodesic is radial near $p$: $q(t)=\exp_p(sv)$ for
a fixed unit vector $v\in T_pM$.

By \eqref{eq:normalized-potential-form} and
Proposition~\ref{prop:isotropy-rigidity}~(i),
\begin{eqnarray*}
 g_p(Du_i,Du_j)=Q(u_i,u_j)=\delta_{ij}.
\end{eqnarray*}
Thus $d\mathbf u_p:T_pM\to\R^n$ is an isometry, and it is
equivariant by Proposition~\ref{prop:isotropy-rigidity}\textup{(iii)}.
Use this identification to choose orthonormal normal coordinates.
Equivariance of $\exp_p$ then gives
$G\cdot q(t)=\exp_p(s\Sph^{n-1})$.
Gauss's lemma and invariance under the stabilizer
give, along the radial axis,
\begin{equation}\label{eq:invariant-center-metric}
 g_{se_1}=dx_1^2+b(s)\sum_{i=2}^n dx_i^2,
\end{equation}
where
\begin{equation}\label{eq:center-radial-normalization}
 b(s)=\frac{\tr_\delta g_{se_1}-1}{n-1},\qquad b(0)=1.
\end{equation}
Here $\delta$ is the Euclidean metric in the normal coordinates.
The function $b$ is smooth and positive near zero.
The angular unit vector $e_2\in T_{e_1}\Sph^{n-1}$ corresponds
to $se_2$ at $se_1$, so
\begin{eqnarray*}
 \rho_*(s)^2=g_{se_1}(se_2,se_2)=s^2b(s).
\end{eqnarray*}
Hence $\rho_*(s)=s\sqrt{b(s)}$ for $s\ge0$, which proves
\textup{(ii)} since $b(0)=1$.
\end{proof}

\section{Classification with compact boundary}
\label{sec:compact-boundary-classification}

The preceding section supplies the global coordinates for both
classifications. We now use the constant sectional curvature already
established in Theorem~\ref{thm:interior-maximality} to determine their
warping functions.

\begin{proposition}\label{prop:warping-equations}
Under the standing assumptions of Section~\ref{sec:orbit-space}, the
functions in \eqref{eq:global-warped-and-potential} satisfy $a=\rho$
and
\begin{equation}\label{eq:space-form-warping-equations}
 \rho''+\kappa\rho=0,\qquad
 (\rho')^2+\kappa\rho^2=1,
 \qquad \kappa=\tfrac{R}{n(n-1)}.
\end{equation}
\end{proposition}

\begin{proof}
The boundary coordinate normalization gives $Q(u_i,u_j)=\delta_{ij}$.
At $q(t)$, where $\mathbf u=a e_1$, this implies
\begin{eqnarray*}
 g(Du_j,Du_k)=\delta_{jk},\qquad X_{1j}=aDu_j
 \qquad(2\le j,k\le n).
\end{eqnarray*}
The fields $X_{1j}$ correspond to the standard unit tangent vectors
$e_j$ at $e_1\in\Sph^{n-1}$. Comparing their lengths with the orbit
metric in \eqref{eq:global-warped-and-potential} gives $a=\rho$.

For this spherical warped product, the radial and tangential sectional
curvatures are \cite[Chapter~7]{ONeill1983}
\begin{eqnarray*}
 K_{\mathrm{rad}}=-\tfrac{\rho''}{\rho},\qquad
 K_{\mathrm{tan}}=\tfrac{1-(\rho')^2}{\rho^2}.
\end{eqnarray*}
Both equal $\kappa$, proving \eqref{eq:space-form-warping-equations}.
\end{proof}

Define
\begin{equation}\label{eq:model-sine}
 s_\kappa(r)=
 \begin{cases}
  \kappa^{-1/2}\sin(\sqrt\kappa\,r),&\kappa>0,\\
  r,&\kappa=0,\\
  (-\kappa)^{-1/2}\sinh(\sqrt{-\kappa}\,r),&\kappa<0.
 \end{cases}
\end{equation}
For $\kappa>0$, the polar coordinate is restricted to
$0\le r<\pi/\sqrt\kappa$.

\begin{proof}[\textbf{\textit{Proof of Theorem~\ref{thm:main}}}]
Suppose that $M$ is compact and $\dim\ker\Phi^*=n$.
By Lemma~\ref{lem:orbit-interval}, its orbit space is $[0,L]$, with
one endpoint corresponding to the selected boundary component.
The other endpoint is a fixed point or a second boundary component.

In the first case, measure arclength $r$ from the fixed point $p$ and
write $f(r)$ for the orbit radius. Lemma~\ref{lem:endpoint-regularity}
and Proposition~\ref{prop:warping-equations} give
\begin{equation}\label{eq:ball-radius-initial-system}
 f''+\kappa f=0,\qquad f(0)=0,\qquad f'(0)=1.
\end{equation}
Hence $f=s_\kappa$ by uniqueness for this initial value problem, and
\begin{eqnarray*}
 g=dr^2+s_\kappa(r)^2g_{\Sph^{n-1}},\qquad 0\le r\le r_+.
\end{eqnarray*}
Positivity of $f$ on $(0,r_+]$ gives $r_+<\pi/\sqrt\kappa$ when
$\kappa>0$. The $G$-translates of the reversed section $q$ are the
radial geodesics from $p$, and they parametrize
$M\setminus\{p\}$ bijectively by Lemma~\ref{lem:global-section}.
The normal coordinates at $p$ and the boundary collar extend this
parametrization smoothly to the center and boundary. Choose
$o\in\Mspace_\kappa^n$ and a linear isometry
$I:T_o\Mspace_\kappa^n\to T_pM$. The matching polar metrics give
the global isometry
\begin{equation}\label{eq:ball-model-isometry}
 \exp_p\circ I\circ\exp_o^{-1}:
 \overline B_{r_+}(o)\longrightarrow M.
\end{equation}

In the second case, $\rho>0$ on $[0,L]$. Fix $t_0\in(0,L)$.
The first integral in \eqref{eq:space-form-warping-equations} gives
$r_*>0$ and $\varepsilon\in\{-1,1\}$ such that
\begin{eqnarray*}
 (\rho(t_0),\rho'(t_0))
 =\bigl(s_\kappa(r_*),\varepsilon s_\kappa'(r_*)\bigr),
\end{eqnarray*}
with $r_*<\pi/\sqrt\kappa$ if $\kappa>0$.
ODE uniqueness gives $\rho(t)=s_\kappa(r_*+\varepsilon(t-t_0))$.
Since $\rho$ stays positive, the image interval $[r_-,r_+]$ lies in
$(0,\infty)$, and in $(0,\pi/\sqrt\kappa)$ when $\kappa>0$.
Changing from $t$ to $r=r_*+\varepsilon(t-t_0)$ therefore gives
\begin{eqnarray*}
 g=dr^2+s_\kappa(r)^2g_{\Sph^{n-1}},\qquad r_-\le r\le r_+.
\end{eqnarray*}
The product map of Lemma~\ref{lem:global-section} includes both
boundary orbits. Comparing it with the polar coordinates of
$\Mspace_\kappa^n$ gives an isometry from the closed concentric
annulus onto $M$.

For the converse, on either model set
$u_i=s_\kappa(r)\theta_i$, $1\le i\le n$.
These are the appropriate ambient coordinate functions in the
Euclidean, spherical, and hyperboloid models. They satisfy
$\Hess u_i=-\kappa u_i g$ and hence $L^*u_i=0$.
On a ball they are smooth at the center; equivalently,
$s_\kappa(r)/r$ is smooth as a function of $r^2$.
On any boundary sphere, writing $\nu=\varepsilon\partial_r$ gives
\begin{eqnarray*}
 h=\varepsilon\tfrac{s_\kappa'}{s_\kappa}\hat g,
 \qquad (u_i)_\nu\hat g
 =\varepsilon s_\kappa'\theta_i\hat g=u_i h.
\end{eqnarray*}
Thus the $u_i$ are b-static potentials. Their restrictions
to a boundary sphere are independent, so the dimension bound in
Proposition~\ref{prop:one-jet} gives $\dim\ker\Phi^*=n$.
\end{proof}

\begin{proof}[\textbf{\textit{Proof of Theorem~\ref{thm:noncompact-main}}}]
Suppose that $M$ is noncompact and $\dim\ker\Phi^*=n$.
All the results of Section~\ref{sec:orbit-space} apply to the selected
compact boundary component. In particular,
\begin{equation}\label{eq:noncompact-orbit-ray}
 M/G=[0,\infty),\qquad\partial M=\Sigma,
\end{equation}
there are no fixed points, and the global product metric is
\begin{equation}\label{eq:noncompact-warped-product}
 g=dt^2+\rho(t)^2g_{\Sph^{n-1}},\qquad \rho(0)>0,
\end{equation}
where $\rho$ is positive on the whole ray and satisfies
\eqref{eq:space-form-warping-equations}.

If $\kappa>0$, every nonzero solution of $\rho''+\kappa\rho=0$
has a zero at positive time, contradicting this positivity.
Hence $\kappa\le0$, and $(\rho')^2=1-\kappa\rho^2\ge1$.
The derivative has constant sign. If it were negative, then
$\rho(t)\le\rho(0)-t$, again a contradiction. Thus $\rho'>0$.
There is a unique $r_0>0$ with $s_\kappa(r_0)=\rho(0)$, and
the first integral gives $\rho'(0)=s_\kappa'(r_0)$.
ODE uniqueness yields
\begin{eqnarray*}
 \rho(t)=s_\kappa(t+r_0)\qquad(t\ge0).
\end{eqnarray*}
Setting $r=t+r_0$ identifies $M$ with
$\Mspace_\kappa^n\setminus B_{r_0}(o)$, including its boundary.

Conversely, this exterior is complete, and the model verification in
the proof of Theorem~\ref{thm:main}, with outward normal
$-\partial_r$ at $r=r_0$, gives $n$ independent b-static
potentials. Proposition~\ref{prop:one-jet} gives equality.
Finally, the compact exhaustion
$[r_0,R]\times\Sph^{n-1}$ has connected complement for every
$R>r_0$, so the exterior has exactly one end.
\end{proof}

\section{Classification with noncompact boundary}
\label{sec:noncompact-boundary-classification}

We now classify complete manifolds of dimension $n\ge3$ with
$\dim\ker\Phi^*=n$ and a noncompact boundary component.
By Corollary~\ref{cor:complete-noncompact-boundary}, all boundary
components are noncompact and $\kappa\le0$.
The ambient curvature and the inertia of $Q$ distinguish
the flat, horospherical, and hyperbolic-boundary cases.
In each case, the potentials define an equivariant map
to which Lemma~\ref{lem:equivariant-product} applies,
giving a global product with a transverse curve.

Completeness determines the possible transverse intervals.
For later use, we record the following standard fact;
cf.~\cite[Problems~13-18 and~15-14]{Lee2013}.

\begin{lemma}\label{lem:complete-transverse-curve}
Let $(M,g)$ be a complete Riemannian manifold with smooth boundary,
and let $C\subset M$ be a closed, connected, embedded
one-dimensional submanifold meeting $\partial M$ transversely,
with $\partial C=C\cap\partial M\ne\varnothing$.
Then $C$ is complete in the induced metric and is isometric to a
closed bounded interval or a closed half-line.
\end{lemma}

\begin{proof}[\textbf{\textit{Proof of Theorem~\ref{thm:noncompact-boundary-main}}}]
Assume $\dim\ker\Phi^*=n$ and put $m=n-1$.
By Theorem~\ref{thm:local-main} and
Corollary~\ref{cor:complete-noncompact-boundary}, $\kappa\le0$ and every
boundary component is noncompact.

We first describe the application of
Lemma~\ref{lem:equivariant-product} in the two flat-boundary cases.
Suppose that an isometric $\R^m$-action $T_s$ preserves
$\partial M$ and that a smooth map $x:M\to\R^m$ satisfies
$x(T_s p)=x(p)+s$. The action is free, so its stabilizers agree with
those of the translation action on $\R^m$. By
Lemma~\ref{lem:equivariant-product}, there is a diffeomorphism
$F:C\times\R^m\to M$, where $C=x^{-1}(0)$, with
\begin{equation}\label{eq:translation-product}
 F(c,s)=T_s c,\qquad
 F^{-1}(p)=\bigl(T_{-x(p)}p,x(p)\bigr).
\end{equation}
Since $M\cong C\times\R^m$ is connected with nonempty
boundary, $C$ is connected and $\partial C\ne\varnothing$.
By Lemma~\ref{lem:complete-transverse-curve}, $C$ is complete and its
arclength parameter ranges over a closed bounded interval or half-line $I$.

If the boundary components are flat,
Lemma~\ref{lem:potential-form-boundary} shows that $Q$ is positive
semidefinite of rank $m$.

\smallskip
\noindent\textit{Case $(i)$: $\kappa=0$.}
By Corollary~\ref{cor:complete-noncompact-boundary},
$k_\Sigma=\lambda_\Sigma=0$ on every boundary component, so each
component is totally geodesic. Constants belong to $\ker\Phi^*$, and the radical
of $Q$ is spanned by $1$. Choose $x_1,\ldots,x_m\in\ker\Phi^*$ with
\[
 Q(x_i,x_j)=\delta_{ij}.
\]
The fields $Y_i=Dx_i$ are parallel, commute, and satisfy
\[
 \inner{Y_i}{Y_j}=\delta_{ij},\qquad Y_i(x_j)=\delta_{ij}.
\]
They are tangent to the whole boundary by the Robin condition and
complete by Lemma~\ref{lem:complete-Killing-integration}. Their flows
thus define the required isometric translation action. In
\eqref{eq:translation-product}, the fields $Y_i$ are an orthonormal
frame for the orbits and are orthogonal to $C$. Taking arclength $t$
on $C$ gives
\[
 g=dt^2+g_{\R^m}\quad\text{on }I\times\R^m.
\]

\smallskip
\noindent\textit{Case $(ii)$: $\kappa<0$ and $k_\Sigma=0$.}
Write $\kappa=-a^2$, where $a>0$. Replace $g$ by $a^2g$, which preserves both kernels and
completeness and normalizes the sectional curvature to $-1$.
Thus $\Hess u=ug$ and $Q(u,v)=\inner{Du}{Dv}-uv$.
Choose $0\ne u_0$ in the radical of $Q$ and $u_1,\ldots,u_m$ so that
\begin{equation}\label{eq:horospherical-potential-normalization}
 Q(u_0,u_i)=0\quad(0\le i\le m),\qquad
 Q(u_i,u_j)=\delta_{ij}\quad(1\le i,j\le m).
\end{equation}
The identity $|Du_0|^2=u_0^2$ and one-jet uniqueness show that $u_0$
has no zero. Choose its sign so that $u_0>0$, and set
\[
 t=\log u_0,\qquad x_i=\tfrac{u_i}{u_0},\qquad
 Y_i=X_{u_0,u_i}=u_0Du_i-u_iDu_0.
\]
The normalization gives
\[
 |Du_0|^2=u_0^2,\qquad
 \inner{Du_0}{Du_i}=u_0u_i,\qquad
 \inner{Du_i}{Du_j}=\delta_{ij}+u_iu_j.
\]
Differentiating $t$ and $x_i$, and expanding the inner products gives
\begin{equation}\label{eq:horospherical-translation-identities}
 \begin{aligned}
 &|Dt|=1,\qquad Y_i(t)=0,\qquad Y_i(x_j)=\delta_{ij},\\
 &\inner{Y_i}{Y_j}=u_0^2\delta_{ij},\qquad Dx_i=u_0^{-2}Y_i.
 \end{aligned}
\end{equation}
By \eqref{eq:pair-bracket-bilinear-form} and
\eqref{eq:horospherical-potential-normalization}, $[Y_i,Y_j]=0$.
These boundary-tangent Killing fields are complete, so their flows
define an isometric $\R^m$-action satisfying
\[
 x(T_s p)=x(p)+s,\qquad t(T_s p)=t(p).
\]
Apply \eqref{eq:translation-product} with $C=x^{-1}(0)$.
By \eqref{eq:horospherical-translation-identities}, $Dt$ is a unit
tangent field along $C$. In an oriented arclength coordinate $\tau$
on $C$, the derivative $dt/d\tau$ takes values in $\{-1,1\}$ and
is continuous, hence constant. Thus $t|_C$ is itself a global
arclength coordinate, with image a closed finite interval or a
closed half-line $I$. Since the action is isometric and preserves
$t$, the same identities give
\[
 g=dt^2+e^{2t}g_{\R^m}\quad\text{on }I\times\R^m.
\]

\smallskip
\noindent\textit{Case $(iii)$: $k_\Sigma<0$.}
Normalize the sectional curvature to $\kappa=-1$.

\smallskip
\noindent\textit{Step 1. The normalized potential map.}
By Lemma~\ref{lem:potential-form-boundary}, the form $Q$ is
nondegenerate with one negative eigenvalue. Choose a basis $u_0,u_1,\ldots,u_m$
of $\ker\Phi^*$ with
\begin{equation}\label{eq:hyperbolic-potential-normalization}
 Q(u_0,u_0)=-1,\quad Q(u_0,u_i)=0,\quad
 Q(u_i,u_j)=\delta_{ij}\quad(1\le i,j\le m).
\end{equation}
Since $|Du_0|^2=u_0^2-1$, the function $u_0$ has no zero; choose its
sign so that $u_0>0$. Write
\[
 \eta(z,z)=-z_0^2+z_1^2+\cdots+z_m^2,\qquad
 U=(u_0,u_1,\ldots,u_m).
\]
The Gram matrix of the gradients is
\[
 \bigl(\inner{Du_i}{Du_j}\bigr)_{0\le i,j\le m}
 =\operatorname{diag}(-1,1,\ldots,1)+UU^{\mathsf T}.
\]
Its positive semidefiniteness gives
$0\le-(1+\eta(U,U))$, so $\eta(U,U)\le-1$. Hence
\begin{equation}\label{eq:normalized-hyperbolic-potential-map}
 \Theta=\frac{U}{\sqrt{-\eta(U,U)}}:
 M\longrightarrow\mathbb H^m_{-1}
 =\{z:\eta(z,z)=-1,\ z_0>0\}
\end{equation}
is well defined.

\smallskip
\noindent\textit{Step 2. The isometric group action.}
Let $G=\mathrm{SO}_0(m,1)$ denote the identity component of
the group of linear transformations preserving
$\eta$, and write $\mathfrak{so}(m,1)$ for its Lie algebra.
Its standard linear action preserves $\mathbb H^m_{-1}$.
By \eqref{eq:pair-action-bilinear-form} and
\eqref{eq:pair-bracket-bilinear-form}, the pair fields define a
Lie algebra homomorphism $A\mapsto X_A$ from
$\mathfrak{so}(m,1)$ to the boundary-tangent Killing fields, with
\begin{equation}\label{eq:hyperbolic-infinitesimal-action}
 dU(X_A)=-AU
 \qquad\text{for every }A\in\mathfrak{so}(m,1),
\end{equation}
in accordance with \eqref{eq:rotation-fundamental-convention}.

Let $\varpi:\widetilde G\to G$ be the simply connected covering
group. By Lemma~\ref{lem:complete-Killing-integration}, the
fields $X_A$ integrate to an isometric $\widetilde G$-action
preserving $\partial M$. Along $p(t)=\exp(-tA)\cdot p$,
\eqref{eq:hyperbolic-infinitesimal-action} gives
\begin{eqnarray*}
 \frac{d}{dt}U(p(t))=-AU(p(t)),\qquad U(p(0))=U(p).
\end{eqnarray*}
Hence $U(\exp(-tA)\cdot p)=e^{-tA}U(p)$.
Since exponentials generate $\widetilde G$ and
$\varpi(\widetilde a)$ preserves $\eta$, for every
$\widetilde a\in\widetilde G$ and $p\in M$ we obtain
\begin{equation}\label{eq:hyperbolic-cover-equivariance}
 \begin{aligned}
 U(\widetilde a\cdot p)
   &=\varpi(\widetilde a)U(p),\\
 \Theta(\widetilde a\cdot p)
   &=\varpi(\widetilde a)\Theta(p).
 \end{aligned}
\end{equation}

To descend the action to $G$, fix a boundary component $\Sigma$ of
curvature $k_\Sigma=-b^{-2}$. By Theorem~\ref{thm:local-main},
restriction identifies $\ker\Phi^*$ with the full Hessian solution
space on $\mathbb H^m_{-b^{-2}}$. In the hyperboloid model of radius
$b$, the ambient coordinate functions form a $Q$-orthonormal basis.
Thus $U|_\Sigma$ differs from this coordinate embedding by a Lorentz
transformation preserving the future cone, and
$\Theta|_\Sigma:\Sigma\to\mathbb H^m_{-1}$ is a diffeomorphism.
By \eqref{eq:hyperbolic-cover-equivariance}, every element of
$\ker\varpi$ fixes $\Sigma$ pointwise and hence acts trivially
on $M$, by the normal-collar argument in
Proposition~\ref{prop:rotation-action}.
The action descends to $G$, and for every $A\in G$,
\begin{equation}\label{eq:hyperbolic-equivariance}
 U(A\cdot p)=AU(p),\qquad
 \Theta(A\cdot p)=A\Theta(p).
\end{equation}

\smallskip
\noindent\textit{Step 3. Stabilizers and the global product.}
At $p\in M$, choose
$u\in\ker\Phi^*$ with $u(p)=1$, using $u_0>0$. The potentials $v$
vanishing at $p$ form an $m$-dimensional space, and their gradients
are linearly independent by one-jet uniqueness. Since
$X_{u,v}(p)=Dv(p)$ and every pair field at $p$ belongs to this
gradient space, every orbit has dimension $m$. By \eqref{eq:hyperbolic-equivariance},
\[
 G_p\subseteq G_{\Theta(p)}\cong\mathrm{SO}(m).
\]
Both groups have dimension $m(m-1)/2$. The closed subgroup $G_p$
is therefore open in $G_{\Theta(p)}$, which is connected since
$m\ge2$. Consequently, $G_p=G_{\Theta(p)}.$

Put $e_0=(1,0,\ldots,0)$, $K=G_{e_0}=\mathrm{SO}(m)$, and
$C=\Theta^{-1}(e_0)$. Identifying $G/K$ with $\mathbb H^m_{-1}$,
Lemma~\ref{lem:equivariant-product} gives the diffeomorphism
\begin{equation}\label{eq:hyperbolic-global-product}
 F:C\times\mathbb H^m_{-1}\longrightarrow M,
 \qquad F(c,Ae_0)=A\cdot c.
\end{equation}
Here $K$ fixes $C$ pointwise and $G$ acts on the second
factor. Since $M\cong C\times\mathbb H^m_{-1}$ is connected with nonempty boundary, $C$ is connected and $\partial C\ne\varnothing$.
By Lemma~\ref{lem:complete-transverse-curve}, $C$ is complete and its
arclength parameter ranges over a closed bounded interval or half-line $I$.

\smallskip
\noindent \textit{Step 4. The warping function.} Along $C$, we have $u_i=0$ for $i\ge1$ and $u_0=s>0$. By
\eqref{eq:hyperbolic-potential-normalization},
\[
 \inner{Du_i}{Du_j}=\delta_{ij},\qquad
 \inner{Du_0}{Du_i}=0\quad(i,j\ge1).
\]
The vectors $Du_i$, $i\ge1$, span the orbit tangent space and are
orthogonal to $C$. At a point of $C$, the boost field
$X_{u_0,u_i}=sDu_i$ maps under $d\Theta$ to the unit vector
$e_i\in T_{e_0}\mathbb H^m_{-1}$. Since $G$ acts isometrically,
\eqref{eq:hyperbolic-global-product} gives
\begin{equation}\label{eq:hyperbolic-cosh-product-metric}
 g=dt^2+s(t)^2g_{\mathbb H^m_{-1}}.
\end{equation}
In particular, the curves with constant second coordinate are unit
geodesics. Along $C$, the static equation and
$|Du_0|^2=u_0^2-1$ give
\[
 s''=s,\qquad (s')^2=s^2-1.
\]
The positive solutions are $s(t)=\cosh(t-t_0)$. Translating the
arclength coordinate and the interval $I$ gives
\[
 g=dt^2+\cosh^2t\,g_{\mathbb H^m_{-1}}.
\]

 Restoring the scale from curvature $-1$ to $\kappa=-a^2$ in the
two negative-curvature cases yields
\[
 dt^2+e^{2at}g_{\R^m},\qquad
 dt^2+\cosh^2(at)g_{\mathbb H^m_{-a^2}}.
\]
In the exponential case, a constant rescaling of the Euclidean
coordinates absorbs the factor $a^{-2}$. The interval remains a
closed finite interval or a closed half-line. This proves the
forward classification.

For the converse, each product is complete, since its
base and fiber are complete and its positive warping factor is
bounded below on compact subintervals of $I$.
The boundary consists of the noncompact endpoint slices.
In the flat case, $1,x_1,\ldots,x_m$ are independent b-static
potentials. In the normalized exponential case, the functions
\[
 e^t, e^t x_1,\ldots,e^t x_m
\]
satisfy $\Hess u=ug$. At an endpoint with outward unit normal
$\nu=\varepsilon\partial_t$, $\varepsilon\in\{-1,1\}$, we have
$h=\varepsilon\hat g$ and $u_\nu=\varepsilon u$, so they also
satisfy the boundary equation. In the normalized hyperbolic-fiber
case, let $y_0,\ldots,y_m$ be the ambient coordinate functions on
$\mathbb H^m_{-1}$. The independent functions
\[
 u_i(t,y)=\cosh t\,y_i\qquad(0\le i\le m)
\]
satisfy $\Hess u_i=u_i g$. At an endpoint,
$h=\varepsilon\tanh t\,\hat g$ and
$(u_i)_\nu=\varepsilon\tanh t\,u_i$, which gives the boundary
equation. Constant rescaling of the metric preserves these
conclusions. Each model therefore has at least $n=m+1$ independent
b-static potentials, and Proposition~\ref{prop:one-jet}
gives equality.
\end{proof}

\section{Classification of surfaces}
\label{sec:surface-classification}

By Remark~\ref{rem:interior-dimension-two}, a single static
potential on a surface already forces constant Gaussian curvature.
The global classification nevertheless requires a separate analysis
of the identifications preserving both potentials.
We carry out this analysis without assuming orientability.
Besides angular translations, the spherical case allows
compositions with a radial reflection.

We first record the boundary equation needed below.
The constancy of $H$ on each connected boundary component
follows from \cite{ShengLocalized2024}.
The boundary Hessian calculation and the one-jet dimension
argument in the proof of Lemma~\ref{lem:boundary-Hessian-system}
also apply in dimension two.
Since $\lambda=H$ and
$\Hess_{\hat g}\phi=\phi''\hat g$ in a local arclength
coordinate, these arguments give the following lemma.

\begin{lemma}\label{lem:surface-local}
Let $(M^2,g)$ be a connected Riemannian surface with nonempty
smooth boundary and $\dim\ker\Phi^*=2$.
Then each boundary component $\Sigma$ has constant geodesic
curvature $H$, and restriction defines an isomorphism from
$\ker\Phi^*$ onto the space of all $\phi\in C^\infty(\Sigma)$
satisfying
\begin{equation}\label{eq:surface-boundary-ode}
 \phi''(s)+(\kappa+H^2)\phi(s)=0,
\end{equation}
where $\kappa=R/2$ and $s$ is a local arclength coordinate
on $\Sigma$.
\end{lemma}

\begin{lemma}\label{lem:surface-product}
Let $(N,g)$ be a complete connected Riemannian surface with nonempty smooth
boundary. Let $X$ be a nowhere vanishing Killing field tangent to
$\partial N$, and let $z\in C^\infty(N)$ satisfy
\begin{equation}\label{eq:surface-flow-coordinate}
 Dz=\frac{X}{|X|^2}.
\end{equation}
Then $Xz=1$, the flow $\varphi_\tau$ of $X$ is complete, and
$C=\{z=0\}$ is a connected closed embedded curve with
$\partial C=C\cap\partial N\ne\varnothing$, isometric to a closed bounded
interval or a closed half-line. Moreover
\begin{equation}\label{eq:surface-global-product}
 F:C\times\mathbb R\longrightarrow N,\qquad F(c,\tau)=\varphi_\tau(c),
\end{equation}
is a diffeomorphism with $z(F(c,\tau))=\tau$, and
\begin{equation}\label{eq:surface-warped-product}
 F^*g=ds^2+f(s)^2d\tau^2,\qquad f=|X|\big|_C>0,
\end{equation}
where $s$ is the arclength coordinate on $C$.
\end{lemma}

\begin{proof}
We have $Xz=\inner{Dz}{X}=1$. By
Lemma~\ref{lem:complete-Killing-integration}\textup{(i)}, $\varphi_\tau$ is
a complete flow of isometries preserving $\partial N$, and
$z(\varphi_\tau p)=z(p)+\tau$. Thus $\tau\cdot p=\varphi_\tau(p)$ is an
isometric $\mathbb R$-action preserving $\partial N$, the map $z$ is
equivariant for translation, and both actions have trivial stabilizers.
Lemma~\ref{lem:equivariant-product} with $G=\mathbb R$ and $K=\{0\}$ shows
that $C$ is a closed embedded curve transverse to $\partial N$, that
$\partial C=C\cap\partial N$, and that $F$ is a diffeomorphism with inverse
$p\mapsto(\varphi_{-z(p)}(p),z(p))$.

The curve $C$ is connected because $N\cong C\times\mathbb R$ is. It meets
the boundary: if $p\in\partial N$, then
$\varphi_{-z(p)}(p)\in C\cap\partial N$. By
Lemma~\ref{lem:complete-transverse-curve}, $C$ is complete and isometric to
a closed bounded interval or a closed half-line.

By \eqref{eq:surface-flow-coordinate}, $dz=\inner{X}{\cdot}/|X|^2$, so
$TC=\ker dz=X^\perp$. Let $c(s)$ be a unit-speed parametrization of $C$.
Then $\partial_sF=d\varphi_\tau(c'(s))$ and $\partial_\tau F=X\circ F$.
Since $\varphi_\tau$ is an isometry and preserves $X$, we obtain
\[
 \inner{\partial_sF}{\partial_\tau F}=\inner{c'(s)}{X(c(s))}=0,
 \qquad |\partial_sF|=1,\qquad |\partial_\tau F|=|X(c(s))|=f(s).
\]
This proves \eqref{eq:surface-warped-product}.
\end{proof}

Recall $s_\kappa$ from \eqref{eq:model-sine}, and put
\begin{equation}\label{eq:surface-polar-functions}
 r_\kappa=
 \begin{cases}
  \pi/\sqrt\kappa,&\kappa>0,\\
  +\infty,&\kappa\le0.
 \end{cases}
\end{equation}

\begin{theorem}\label{thm:surface-classification}
Let $(M^2,g)$ be a complete connected Riemannian surface with
nonempty smooth boundary. Then $\dim\ker\Phi^*=2$ if and only if
$(M,g)$ is isometric to one of the following models.
\begin{enumerate}[label=\textup{(\roman*)}]
\item\label{item:surface-disks}
A closed geodesic disk of radius $b$, $0<b<r_\kappa$, in
$\mathbb M_\kappa^2$.

\item\label{item:surface-polar}
A warped product
\begin{equation}\label{eq:surface-polar-model}
 I\times\Theta,\qquad
 g=dr^2+s_\kappa(r)^2d\theta^2,
\end{equation}
where $\Theta=\mathbb R$ or
$\Theta=\mathbb R/(2\pi m\mathbb Z)$ for an integer $m\ge1$, and
\begin{eqnarray*}
 I=[r_-,r_+],\quad 0<r_-<r_+<r_\kappa,
\end{eqnarray*}
or, when $\kappa\le0$,
\begin{eqnarray*}
 I=[r_-,\infty),\quad r_->0.
\end{eqnarray*}
These are the connected Riemannian covers of concentric annuli
and, for $\kappa\le0$, of geodesic-ball exteriors.

\item\label{item:surface-flat}
A flat product
\begin{equation}\label{eq:surface-flat-model}
 I\times\mathbb R,\qquad g=dt^2+dx^2,
\end{equation}
where $I$ is a nondegenerate closed bounded interval or a closed
half-line.

\item\label{item:surface-hyperbolic}
A hyperbolic warped product of one of the forms
\begin{equation}\label{eq:surface-hyperbolic-models}
 \begin{aligned}
 I\times\mathbb R,&\qquad g=dt^2+e^{2at}dx^2,\\
 I\times\mathbb R,&\qquad g=dt^2+\cosh^2(at)dx^2,
 \end{aligned}
\end{equation}
where $a>0$, $\kappa=-a^2$, and $I$ is a
nondegenerate closed bounded interval or $[A,\infty)$, $A\in\R$.
For the exponential model, $I=(-\infty,B]$, $B\in\R$, is also allowed.

\item\label{item:surface-moebius}
A spherical M\"obius band of curvature $\kappa>0$,
given by the quotient of
\begin{equation}\label{eq:surface-moebius-cover}
 [-b,b]\times\mathbb R,
 \qquad g=dt^2+\frac{\cos^2(\sqrt\kappa\,t)}{\kappa}\,d\theta^2,
 \qquad 0<b<\frac{\pi}{2\sqrt\kappa},
\end{equation}
under the infinite cyclic group generated by
\begin{equation}\label{eq:surface-moebius-generator}
 (t,\theta)\longmapsto(-t,\theta+2\pi m),
 \qquad m\in\mathbb N,\quad m\ge1.
\end{equation}
\end{enumerate}
The nonorientable models are precisely those in
\textup{(v)}.
\end{theorem}

We begin the proof of Theorem~\ref{thm:surface-classification} with the case in which $Q$ is positive definite.
\begin{lemma}\label{lem:surface-positive-classification}
Let $(M^2,g)$ be a complete connected Riemannian surface with nonempty
smooth boundary and $\dim\ker\Phi^*=2$. If the bilinear form $Q$ in
\eqref{eq:potential-bilinear-form} is positive definite on $\ker\Phi^*$,
then $(M,g)$ is isometric to one of the models \textup{(i)},
\textup{(ii)}, or \textup{(v)} in
Theorem~\ref{thm:surface-classification}.
\end{lemma}

\begin{proof}
By Remark~\ref{rem:interior-dimension-two}, $g$ has constant
Gaussian curvature $\kappa$. Choose a $Q$-orthonormal basis $u,v$ of $\ker\Phi^*$ and set
$X=uDv-vDu$. Proposition~\ref{prop:pair-Killing} shows that $X$ is
Killing and tangent to $\partial M$. The normalization gives
$|Du|^2=1-\kappa u^2$, $|Dv|^2=1-\kappa v^2$, and
$\inner{Du}{Dv}=-\kappa uv$. Together with
\eqref{eq:pair-action-bilinear-form}, these identities yield
\begin{equation}\label{eq:surface-positive-identities}
 Xu=-v,\qquad Xv=u,\qquad |X|^2=u^2+v^2 .
\end{equation}

Let $\pi:\widetilde M\to M$ be the universal Riemannian covering, retaining
$u,v,X$ for the lifts. Then $\widetilde M$ is simply connected,
hence orientable, and
$\partial\widetilde M=\pi^{-1}(\partial M)\ne\varnothing$. To see completeness, let $p_j$ be Cauchy in
$\widetilde M$. Completeness of $M$ gives $\pi(p_j)\to q$.
Choose an evenly covered neighborhood $U\supset B_{2\varepsilon}(q)$.
For large $j,k$, $\pi(p_j)\in B_\varepsilon(q)$ and a curve from $p_j$
to $p_k$ of length less than $\varepsilon$ projects into $U$.
Thus the tail lies in one sheet and converges to the lift of $q$ there. The lifts $u,v,X$ are deck-invariant, and
Lemma~\ref{lem:complete-Killing-integration} gives a complete flow
$\varphi_\tau$ of $X$. An isometry is determined by its value and
differential at one point, by the exponential map in the interior and
the boundary collar argument in Proposition~\ref{prop:rotation-action}.

\smallskip\noindent\emph{$(a)$ Common zero.}
Suppose $u(p)=v(p)=0$. Boundary one-jet injectivity implies
$p\notin\partial\widetilde M$, since otherwise both one-jets would lie in
the line $\{0\}\oplus T_p^*(\partial\widetilde M)$. By
\eqref{eq:surface-positive-identities}, the same applies to every zero
of $X$. At any such point, $Du,Dv$ are orthonormal and
\[
 D_ZX=du(Z)Dv-dv(Z)Du.
\]
Thus $d\varphi_\tau|_p$ is rotation through angle $\tau$. Isometry
rigidity gives $\varphi_{2\pi}=\mathrm{id}$, with no smaller positive
period. The resulting effective isometric $\mathbb R/2\pi\mathbb Z$-action
preserves orientation, since each $\varphi_\tau$ is isotopic to the
identity. An isotropy element at $q$ with $X(q)\ne0$ fixes
$X(q)$, so its orientation-preserving differential is the identity;
isometry rigidity makes the element trivial. All nonfixed orbits are
therefore free circles.

As in Lemma~\ref{lem:orbit-interval}, tubular neighborhoods of free
orbits, boundary collars, and exponential disks at fixed points identify
the orbit space as a connected Hausdorff one-dimensional manifold whose
endpoints are precisely the fixed points and boundary orbits. The fixed
point $p$ and a boundary orbit give two distinct endpoints; hence the
orbit space is a compact interval, $p$ is the unique fixed point, and
there is exactly one boundary orbit. The quotient map for this compact
circle action is proper, so
$\widetilde M$ is compact.

For a unit-speed geodesic $\gamma$ issuing from $p$, the Killing equation
gives
\[
 \tfrac{d}{dr}\inner{\dot\gamma}{X}
 =\inner{\dot\gamma}{D_{\dot\gamma}X}=0,
 \qquad \inner{\dot\gamma}{X}(0)=0.
\]
Hence $\gamma$ is orthogonal to the orbits.
Its projection parametrizes the orbit interval monotonically at unit
speed, meeting each regular orbit once and reaching the boundary at the
interval length $b$. Freeness makes
$(r,\theta)\mapsto\varphi_\theta(\gamma(r))$ injective for $0<r\le b$;
the boundary collar and exponential disk give global polar coordinates
\[
 g=dr^2+f(r)^2d\theta^2,\qquad 0\le r\le b,\quad
 \theta\in\mathbb R/2\pi\mathbb Z,\qquad f=|X|.
\]
Smoothness at $p$ and constant curvature imply
\[
 f''+\kappa f=0,\qquad f(0)=0,\qquad f'(0)=1.
\]
Hence $f=s_\kappa$, and $f>0$ on $(0,b]$ gives $b<r_\kappa$.
Every deck transformation fixes the unique zero $p$ of $X$ and is
therefore the identity. Thus $M=\widetilde M$ is model~\textup{(i)}.

\smallskip\noindent\emph{$(b)$ No common zero.}
Set $f=|X|=(u^2+v^2)^{1/2}>0$. Simple connectedness gives a smooth real
lift $\theta$ of $(u+iv)/f$, satisfying
\begin{equation}\label{eq:surface-positive-angle}
 u=f\cos\theta,\quad v=f\sin\theta,\qquad
 D\theta=\frac{uDv-vDu}{f^2}=\frac{X}{|X|^2}.
\end{equation}
Lemma~\ref{lem:surface-product}, with $z=\theta$, yields
$\widetilde M\cong I\times\mathbb R$ and
$g=dr^2+f(r)^2d\theta^2$, where $r$ is arclength on
$C=\{\theta=0\}$. On $C$, we have $u=f$ and $\partial_\theta u=0$.
Thus constant curvature and $Q(u,u)|_C=1$ give
\begin{equation}\label{eq:surface-positive-odes}
 \left\{
 \begin{aligned}
  f''+\kappa f&=0,\\
  (f')^2+\kappa f^2&=1.
 \end{aligned}
 \right.
\end{equation}
After translating and, if necessary, reflecting $r$, we have
$f=s_\kappa(r)$ on $I\subset(0,r_\kappa)$. Completeness of $C$ includes
every finite endpoint in $I$; positivity of $f$ excludes endpoints at
zeros of $s_\kappa$. Consequently, $I=[r_-,r_+]$ with
$0<r_-<r_+<r_\kappa$, or, when $\kappa\le0$, $I=[r_-,\infty)$ with $r_->0$.

To determine the deck group $\Gamma$ of $\pi$,
let $\gamma\in\Gamma$.
Since $\gamma$ preserves $u$ and $v$, connectedness gives
\[
 f\circ\gamma=f,\qquad
 \theta\circ\gamma=\theta+2\pi\ell_\gamma,
 \qquad \ell_\gamma\in\mathbb Z.
\]
Since $\gamma^*g=g$, $f\circ\gamma=f$, and
$d(\theta\circ\gamma)=d\theta$, we obtain
\[
 d(r\circ\gamma)\otimes d(r\circ\gamma)=dr\otimes dr,
 \qquad r\circ\gamma=\varepsilon r+c,
 \qquad \varepsilon\in\{\pm1\},
\]
where $\varepsilon$ and $c$ are constant. Using $f=s_\kappa$ on
$I\subset(0,r_\kappa)$, we conclude that
\begin{equation}\label{eq:surface-radial-deck}
 \begin{aligned}
 \gamma(r,\theta)&=(r,\theta+2\pi\ell_\gamma),\quad\text{or}\\
 \gamma(r,\theta)&=\left(\tfrac{\pi}{\sqrt\kappa}-r,
\theta+2\pi\ell_\gamma\right).
 \end{aligned}
\end{equation}
The second form requires $\kappa>0$ and $I$ symmetric about
$\pi/(2\sqrt\kappa)$.

The homomorphism $\Gamma\to\mathbb Z$, $\gamma\mapsto\ell_\gamma$,
is injective. Indeed, an element in its kernel is the identity or
the pure radial reflection, which fixes the points with
$r=\pi/(2\sqrt\kappa)$; a nontrivial deck transformation has no fixed
point. Thus $\Gamma$ is trivial or infinite cyclic, generated by
$\gamma_1$ with $\ell_{\gamma_1}=m\ge1$. The cases
\[
 M\cong
 \begin{cases}
  I\times\mathbb R,&\Gamma=\{e\},\\
  I\times(\mathbb R/2\pi m\mathbb Z),
    &\Gamma=\langle\gamma_1\rangle,\quad r\circ\gamma_1=r,
 \end{cases}
\]
give \textup{(ii)}. If $\gamma_1$ reverses $r$, then $\kappa>0$ and
\[
 I=[r_0-b,r_0+b],\qquad
 r_0=\tfrac{\pi}{2\sqrt\kappa},\qquad 0<b<r_0.
\]
The coordinate $t=r-r_0\in[-b,b]$ transforms the metric and generator
into \eqref{eq:surface-moebius-cover} and
\eqref{eq:surface-moebius-generator}, respectively, giving
\textup{(v)}.
\end{proof}

\begin{proof}[\textbf{Proof of
Theorem~\ref{thm:surface-classification}}]
Assume $\dim\ker\Phi^*=2$. Remark~\ref{rem:interior-dimension-two}
and Lemma~\ref{lem:surface-local} give constant Gaussian curvature $\kappa$
and constant $H$ on each boundary component. Let $Q$ be the constant bilinear form
\eqref{eq:potential-bilinear-form}. For $\phi=u|_\Sigma$ and
$\psi=v|_\Sigma$, the Robin condition gives
\begin{equation}\label{eq:surface-boundary-form}
 Q(u,v)=\phi'\psi'+(\kappa+H^2)\phi\psi.
\end{equation}
The boundary one-jet isomorphism identifies $Q$ with
\[
 (s,\xi)\longmapsto(\kappa+H^2)s^2+\xi^2.
\]
Thus $Q$ is positive definite, positive semidefinite of rank one,
or of signature $(1,1)$. The positive-definite case follows from
Lemma~\ref{lem:surface-positive-classification}.

For each basis $u,v$ used in the remaining cases, write $X=X_{u,v}$
as in \eqref{eq:intro-Xuv}. By Proposition~\ref{prop:pair-Killing},
$X$ is Killing and tangent to $\partial M$, and its action on
potentials is given by \eqref{eq:pair-action-bilinear-form}.

\medskip\noindent\textbf{Step 1: $Q$ semidefinite of rank one.}
Choose a basis $u,v$ with $Q(u,u)=Q(u,v)=0$ and $Q(v,v)=1$.
By \eqref{eq:surface-boundary-form}, $\kappa=-H^2\le0$, and
\[
 |Du|^2=-\kappa u^2,\qquad
 \inner{Du}{Dv}=-\kappa uv,\qquad
 |Dv|^2=1-\kappa v^2.
\]
If $u(p)=0$, then $Du(p)=0$; one-jet uniqueness would give
$u\equiv0$. Thus $u$ has constant sign. Normalize $u>0$, with
$u\equiv1$ when $\kappa=0$, and set $x=v/u$. Then
\[
 Xu=0,\qquad |X|^2=u^2,\qquad
 Dx=\frac{uDv-vDu}{u^2}=\frac{X}{|X|^2}.
\]
Lemma~\ref{lem:surface-product}, with $z=x$, gives
$M\cong I\times\mathbb R$ and $g=ds^2+u(s)^2dx^2$, where $s$ is
arclength on $C=\{x=0\}$. For $\kappa=0$,
$t=s$ gives \textup{(iii)}. For $\kappa=-a^2<0$, set
$t=a^{-1}\log u$. On $C$,
\[
 X(t)=0,\qquad Dt\in X^\perp=TC,\qquad
 |Dt|^2=\frac{|Du|^2}{a^2u^2}=1.
\]
Hence $t$ is an arclength coordinate on $C$, and
$u=e^{at}$ gives $g=dt^2+e^{2at}dx^2$, the first model in
\textup{(iv)}. Both half-line orientations occur, since
$t\mapsto-t$ does not preserve this metric.

\medskip\noindent\textbf{Step 2: $Q$ of signature $(1,1)$.}
By \eqref{eq:surface-boundary-form}, $\kappa<-H^2\le0$.
Write $\kappa=-a^2$ with $a>0$, and choose a basis with
$Q(u,u)=-1$, $Q(v,v)=1$, and $Q(u,v)=0$. Then
\[
 |Du|^2=a^2u^2-1,\qquad
 \inner{Du}{Dv}=a^2uv,\qquad |Dv|^2=1+a^2v^2.
\]
In particular, $a^2u^2\ge1$, so $u$ has no zero; take $u>0$.
Moreover,
\[
 0\le |Du|^2|Dv|^2-\inner{Du}{Dv}^2=a^2(u^2-v^2)-1,
\]
so $f:=\sqrt{u^2-v^2}\ge a^{-1}$ and $u>|v|$.
Set $y=\operatorname{arctanh}(v/u)$, so $u=f\cosh y$ and
$v=f\sinh y$. Since $Xu=v$ and $Xv=u$,
\[
 |X|^2=f^2,\qquad Xy=1,\qquad Dy=\tfrac{X}{|X|^2}.
\]
Lemma~\ref{lem:surface-product}, with $z=y$, gives
$M\cong I\times\mathbb R$ and $g=dt^2+f(t)^2dy^2$, where $t$ is
arclength on $C=\{y=0\}$. On $C$, we have $u=f$ and
$\partial_yu=0$. Constant curvature and $Q(u,u)|_C=-1$ therefore give
\[
 \left\{
 \begin{aligned}
  f''&=a^2f,\\
  (f')^2-a^2f^2&=-1.
 \end{aligned}
 \right.
\]
Writing $f=Ae^{at}+Be^{-at}$ gives $AB=(4a^2)^{-1}$, and $f>0$
forces $A,B>0$. Translating $t$ and setting $x=y/a$ yields
\[
 f=a^{-1}\cosh(at),\qquad g=dt^2+\cosh^2(at)\,dx^2,
\]
the second model in \textup{(iv)}. The reflection
$(t,x)\mapsto(-t,x)$ is an isometry, so $I=(-\infty,B]$
may be replaced by $[-B,\infty)$.

Lemma~\ref{lem:surface-positive-classification}
and Steps~1--2 establish the forward implication.
Models \textup{(i)}--\textup{(iv)} are orientable disks or products,
whereas the generator \eqref{eq:surface-moebius-generator} reverses
orientation. Thus precisely the models in \textup{(v)} are
nonorientable.

\medskip\noindent\textbf{Step 3: the converse.}
The disks and the M\"obius bands are compact. Every remaining model
is $I\times F$ with
\[
 g=dr^2+f(r)^2dx^2,
 \qquad F=\mathbb R\ \text{or}\ \mathbb R/2\pi m\mathbb Z,
\]
where $I$ is a closed interval or half-line and $f$ is continuous
and positive. Let $p_j=(r_j,x_j)$ be a Cauchy sequence. Since
\[
 |r_j-r_k|\le d(p_j,p_k),
\]
we have $r_j\to r_\infty\in I$. Set
\[
 J=I\cap[r_\infty-1,r_\infty+1],
 \qquad f_{\min}=\min_J f>0.
\]
For large $j,k$, every curve joining $p_j$ to $p_k$ with length
less than $1/2$ lies in $J\times F$. Taking the infimum over such
curves gives
\[
 f_{\min}\,d_F(x_j,x_k)\le d(p_j,p_k).
\]
Thus $x_j$ is Cauchy in the complete fiber $F$, so
$x_j\to x_\infty$ and $p_j\to(r_\infty,x_\infty)$.

For the kernel dimension, take the following pairs of potentials.
\begin{equation}\label{eq:surface-model-potentials}
\setlength{\arraycolsep}{5pt}
\renewcommand{\arraystretch}{1.15}
\begin{array}{@{}lcc@{}}
 \hline
 \text{model} & \multicolumn{2}{c}{\text{potentials}}\\ \hline
 \text{disk or polar}&s_\kappa(r)\cos\theta&s_\kappa(r)\sin\theta\\
 \text{flat product}&1&x\\
 \text{exponential}&e^{at}&e^{at}x\\
 \text{hyperbolic cosine}&\cosh(at)\cosh(ax)&\cosh(at)\sinh(ax)\\
 \text{spherical M\"obius}&\cos(\sqrt\kappa\,t)\cos\theta
   &\cos(\sqrt\kappa\,t)\sin\theta\\ \hline
\end{array}
\end{equation}
The polar pair extends smoothly over a disk center and is
$2\pi m$-periodic in $\theta$. The last pair is invariant under
\eqref{eq:surface-moebius-generator}, since $\cos(\sqrt\kappa\,t)$
is even. Thus all pairs are globally defined.

Away from a disk center, write $g=dt^2+f(t)^2dx^2$, with
$(t,x)=(r,\theta)$ in the polar models and $x=\theta$ in the
M\"obius models. Up to nonzero constant factors, each listed
potential is $w=f(t)\psi(x)$, where
\[
 f''+\kappa f=0,\qquad (f')^2+\kappa f^2=c,\qquad \psi''+c\psi=0.
\]
Here $c=1$ for the polar and M\"obius models, $c=0$ for the flat
and exponential models, and $c=-a^2$ for the hyperbolic cosine model.
Consequently,
\begin{equation}\label{eq:surface-separated-Hessian}
 \begin{aligned}
 \Hess w(\partial_t,\partial_t)&=f''\psi=-\kappa w,\\
 \Hess w(\partial_t,\partial_x)&=f'\psi'-\tfrac{f'}{f}f\psi'=0,\\
 \Hess w(\partial_x,\partial_x)&=f\psi''+f(f')^2\psi=f\bigl((f')^2-c\bigr)\psi=-\kappa f^2w.
 \end{aligned}
\end{equation}
Hence $\Hess w=-\kappa wg$, also at a disk center by smoothness.
At each endpoint,
\[
 \nu=\varepsilon\partial_t,\qquad
 H=\varepsilon\tfrac{f'}{f},\qquad
 w_\nu=\varepsilon f'\psi=Hw,
 \qquad \varepsilon\in\{\pm1\}.
\]
Each pair consists of two linearly independent b-static potentials, so
the dimension bound in
Proposition~\ref{prop:one-jet} yields $\dim\ker\Phi^*=2$.
\end{proof}

The compact models in
Theorem~\ref{thm:surface-classification} are precisely the space-form
geodesic disks, finite cyclic Riemannian covers of concentric
space-form annuli, and spherical M\"obius bands.
Only the M\"obius bands are nonorientable.
Each compact model has at most two boundary components: one for disks and
M\"obius bands, and two for annular covers.

\begin{remark}\label{rem:surface-boundary-modes}
For a compact boundary component $\Sigma$ of length $L$,
\eqref{eq:surface-boundary-ode} has a two-dimensional space of
periodic solutions. Hence $\kappa+H^2>0$ and
\begin{equation}\label{eq:surface-boundary-period}
 \sqrt{\kappa+H^2}\,L=2\pi q,
 \qquad \kappa+H^2=q^2\lambda_1(\Sigma),
\end{equation}
where $q\ge1$ is an integer and $\lambda_1(\Sigma)$ is the first
nonzero eigenvalue of $-\Delta_\Sigma$. The restrictions form the
$q$-th Fourier mode space, not necessarily the first eigenspace.
The mode order is $1$ for a disk, $m$ for an annular cover with
angular period $2\pi m$, and $2m$ for the M\"obius model with
generator $(t,\theta)\mapsto(-t,\theta+2\pi m)$, whose single
boundary has angular period $4\pi m$.
\end{remark}

Several higher-dimensional conclusions fail
for surfaces. Annular covers with $m\ge2$ cannot embed isometrically
in $\mathbb M^2_\kappa$, since their boundaries are longer than
geodesic circles of the same geodesic curvature.
M\"obius bands are excluded by nonorientability.
Moreover, spherical annuli have complete noncompact universal covers
of positive curvature with noncompact boundary, while the flat strip
$[0,1]\times\mathbb R$ has two ends.

\section{Sharpness and the role of completeness}
\label{sec:sharpness}

The dimension assumption in the classification cannot be weakened to
$\dim\ker\Phi^*=n-1$. The first example shows that constant sectional
curvature can fail at this dimension, even for a compact manifold with
connected totally geodesic boundary. The second shows that the number
of boundary components is then unbounded, even for Euclidean domains.

\begin{example}\label{ex:submaximal-product}
Let $n\ge3$, and consider the product
\[
 M=\Sph^2_+\bigl((n-2)^{-1/2}\bigr)
          \times\Sph^{n-2}(1),
 \qquad g=g_1+g_2,
\]
where the first factor is the closed round hemisphere. For $n=3$,
the second factor is the unit circle. Then $M$ is compact and
connected, its boundary is connected and totally geodesic, and
\(
 \dim\ker\Phi^*=n-1,
\)
although $g$ does not have constant sectional curvature.

Indeed, the Ricci and scalar curvatures are
\[
 \Ric=(n-2)g_1+(n-3)g_2,
 \qquad R=(n-1)(n-2),
\]
so
\[
 \Ric-\tfrac{R}{n-1}g=-g_2.
\]
Let $x_1,\ldots,x_{n-1}$ be the ambient coordinate functions on
$\Sph^{n-2}(1)\subset\R^{n-1}$, viewed as functions on $M$.
They satisfy
\[
 \Hess x_i=-x_i g_2
 =x_i\left(\Ric-\tfrac{R}{n-1}g\right),
\]
and hence $L^*x_i=0$. When $n=3$, these are the two functions
$\cos\theta$ and $\sin\theta$ on the circle.
The boundary is the product of the equator of the first factor with
the second factor. Therefore $h=0$ and $H=0$; its outward normal is
tangent to the first factor, so
\[
 (x_i)_\nu=0=\tfrac{H}{n-1}x_i,
 \qquad (x_i)_\nu\hat g=x_i h.
\]
Thus the $n-1$ linearly independent functions $x_i$ belong to
$\ker\Phi^*$. A two-plane tangent to the first factor has sectional
curvature $n-2$, whereas a mixed two-plane has sectional curvature
zero. By Theorem~\ref{thm:interior-maximality}, this excludes
$\dim\ker L^*\ge n$. Consequently,
\[
 \ker\Phi^*=\ker L^*=\Span\{x_1,\ldots,x_{n-1}\}.
\]
This construction restricts familiar closed product examples to a
hemisphere; compare the product models in
\cite[Corollaries~4.2--4.3]{HePetersenWylie2015}.
\end{example}

\begin{example}\label{ex:submaximal-many-components}
Fix $n\ge2$ and an integer $N\ge1$. Choose distinct nonzero numbers
$t_1,\ldots,t_N$, set $c_j=t_je_1\in\R^n$, and choose $R>0$ and
$\varepsilon_j>0$ so that the closed balls
$\overline{B_{\varepsilon_j}(c_j)}$ are pairwise disjoint and contained
in $B_R(0)$. The Euclidean domain
\[
 \Omega=\overline{B_R(0)}
          \setminus\bigcup_{j=1}^N B_{\varepsilon_j}(c_j)
\]
is compact and connected, and its boundary has $N+1$ components.
However,
\[
 \ker\Phi^*=\Span\{x_2,\ldots,x_n\},
 \qquad \dim\ker\Phi^*=n-1.
\]

To verify this, recall that the interior static potentials on a
connected Euclidean domain are precisely the affine functions
$u=a+b\cdot x$. Every boundary component is totally umbilic, so
the tensor boundary condition is equivalent to
$u_\nu=Hu/(n-1)$ on that component. On the outer sphere,
\[
 \nu=\tfrac{x}{R},\qquad \tfrac{H}{n-1}=\tfrac1R,
 \qquad u_\nu-\tfrac{H}{n-1}u=-\tfrac aR.
\]
On the sphere bounding the $j$th removed ball, the normal outward
from $\Omega$ points into that ball. Thus
\[
 \nu=-\tfrac{x-c_j}{\varepsilon_j},
 \qquad \tfrac{H}{n-1}=-\tfrac1{\varepsilon_j},
 \qquad
 u_\nu-\tfrac{H}{n-1}u
 =\tfrac{a+b\cdot c_j}{\varepsilon_j}.
\]
The boundary conditions are therefore
\[
 a=0,\qquad a+t_jb_1=0\quad(1\le j\le N).
\]
Since $N\ge1$ and each $t_j\ne0$, these conditions are equivalent
to $a=b_1=0$, proving the asserted kernel formula.
\end{example}

The following example shows that, without completeness, a manifold
with $\dim\ker\Phi^*=n$ may have arbitrarily many
boundary components.
\begin{example}[Arbitrarily many boundary components]
	\label{ex:incomplete-many-components}
	Let $U\subset\R^{n-1}$ be a nonempty open set and equip
	\begin{eqnarray*}
		\Omega_U=
		\bigl(\R^{n-1}\times(0,\infty)\bigr)
		\cup\bigl(U\times\{0\}\bigr)
	\end{eqnarray*}
	with the Euclidean metric. The set $\Omega_U$ is open in the smooth
	manifold with boundary $\R^{n-1}\times[0,\infty)$, because its
	complement is the closed set
	$(\R^{n-1}\setminus U)\times\{0\}$.
	Thus $\Omega_U$ is a connected smooth manifold with boundary
	$\partial\Omega_U=U\times\{0\}$.
	
	Every interior static potential is affine. Since the boundary is
	totally geodesic and its outward normal is $-\partial_t$, the
	boundary equation reduces to $u_t=0$.
	Consequently, the coefficient of $t$ in an affine potential must
	vanish, and we obtain
	\begin{eqnarray*}
		\ker\Phi^*=\Span\{1,x_1,\ldots,x_{n-1}\},
		\qquad \dim\ker\Phi^*=n.
	\end{eqnarray*}
	The number of connected boundary components is therefore the number
	of connected components of $U$, which can be any positive integer
	or countably infinite. If $U\ne\R^{n-1}$, choose an omitted boundary
	point and approach it along a vertical segment. This produces a
	Cauchy sequence in $\Omega_U$ with no limit in $\Omega_U$, so the
	metric is incomplete. Hence the assumption $\dim\ker\Phi^*=n$
	alone does not imply the ball-or-annulus classification or a bound
	of two boundary components.
\end{example}

\begin{remark}\label{rem:noncompact-exterior-hypotheses}
For $0<a<b<\infty$, the flat manifold
$\{x\in\R^n:a\le|x|<b\}$ has one compact boundary component and
$\dim\ker\Phi^*=n$, but is incomplete at the omitted outer sphere.
Thus a compact boundary component alone does not suffice for
Theorem~\ref{thm:noncompact-main}.
\end{remark}

\end{document}